\documentclass{article}
\usepackage{latexsym,amssymb,amsmath,mathrsfs,a4wide}
\usepackage{hyperref}
\newcommand{\p}{\partial}
\newcommand{\f}{\frac}
\newtheorem{theorem}{Theorem}[section]
\newtheorem{lemma}[theorem]{Lemma}

\newtheorem{definition}{Definition}[section]
\newtheorem{assumption}{Assumption}[section]

\newtheorem{remark}{Remark}[section]

\title{The steady Navier-Stokes and Stokes systems with mixed boundary conditions including one-sided leaks and pressure}
\author{Tujin Kim \\
\thanks{Partially supported by AMSS in Chinese Academy of Sciences.}\\
{\small Institute of Mathematics, }\\
{\small Academy of Sciences,}\\
{\small Pyongyang, DPR Korea } \\
\vspace{.5cm}{\small e-mail: math.tujin@star-co.net.kp}\\
Daomin Cao \\
\thanks{Partially supported by NSFC grants (No. 11271354 and No. 11331010) and Beijing
Center for Mathematics and Information Interdisciplinary Sciences.}\\
{\small Institute of Applied Mathematics, }\\
{\small AMSS, Chinese Academy of Sciences,}\\
{\small Beijing 100190, P. R. China} \\
{\small e-mail: dmcao@amss.ac.cn}
 }

\date{}

\begin{document}

\maketitle

\tableofcontents
\begin{abstract}
In this paper we are concerned with the steady Navier-Stokes and
Stokes problems with mixed boundary conditions involving Tresca
slip, leak condition, one-sided leak conditions, velocity, pressure,
rotation, stress and normal derivative of velocity together. Relying
on the relations among strain, rotation, normal derivative of
velocity and shape of boundary surface, we have variational
formulations for the problems, which consist of five formulae with
five unknowns. We get the variational inequalities equivalent to the
formulated variational problems, which have one unknown. Then, we
study the corresponding variational inequalities and relying the
results for variational inequalities, we get existence, uniqueness
and estimates of solutions to the Navier-Stokes and Stokes problems
with the boundary conditions. Our estimates for solutions
do not depend on the thresholds for slip and leaks.\\
\noindent\textbf{MS Classification 2010:} 35Q30, 35J87, 76D03,
76D05, 49J40
 \vspace{.20cm}

\end{abstract}

\section{Introduction}
\setcounter{equation}{0}

As mathematical models of steady flows of incompressible viscous
Newtonian fluids the Stokes equations
\begin{equation}\label{1.1}
 -\nu\Delta v+\nabla p=f, \,\,\,
 \nabla \cdot v=0\quad \mbox{in}\,\, \Omega,
 \end{equation}
 and Navier-Stokes equations
\begin{equation}\label{1.2}
 -\nu\Delta v+(v\cdot\nabla)v+\nabla p=f, \,\,\,
 \nabla \cdot v=0\quad \mbox{in}\,\, \Omega,
 \end{equation}
are used. For these systems, different natural and artificial
boundary conditions are considered(cf. Introduction of \cite{kc} and
references therein).

Recently several papers are devoted to problems with Tresca slip
boundary condition or leak boundary condition. All these boundary
conditions are called the boundary conditions of friction type,
which are nonlinear.

Tresca slip boundary condition (threshold slip condition) means that
if absolute value of tangent stress on a boundary is less than a
given threshold, then there is not any slip on the boundary surface,
but the absolute value is same with the threshold, then slip on the
boundary surface may occur. Physical and experimental backgrounds of
such boundary conditions are mentioned in several papers(cf.
\cite{f1}, \cite{bb}, \cite{ags}, and especially \cite{hl}). When
$v$ is a solution to \eqref{1.1} or \eqref{1.2}, the strain tensor
is one with the components $\varepsilon_{ij}(v)=\f{1}{2}(\p_{x_i}
v_j+\p_{x_j} v_i)$ and stress tensor $S(v,p)$ is one with components
$S_{ij}=-p\delta_{ij}+2\nu\varepsilon_{ij}(v).$  Let $n$ be the
outward normal unit vector on a boundary surface and $\tau$ tangent
vectors. Then, stress vector on the surface is $\sigma(v,p)=Sn$ and
normal stress $\sigma_n(v,p)=\sigma\cdot n$. Under such notations
Tresca slip boundary condition is expressed by
\begin{equation}\label{1.3}
 |\sigma_\tau(v)|\leq g_\tau,\quad \sigma_\tau(v)\cdot v_\tau+g_\tau
|v_\tau|=0,
 \end{equation}
where and in what follows $\sigma_\tau=\sigma-\sigma_n n$ and
$v_\tau=v-(v\cdot n)n$.

Leak boundary condition means that if absolute value of normal
stress on a boundary is less than a given threshold, then there is
not any leak through the boundary surface, but the absolute value is
same with the threshold, then leak through the boundary surface may
occur. For physical backgrounds of this boundary condition refer to
\cite{f1}, \cite{f4}, \cite{all}. Under notations above leak
boundary condition is expressed by
\begin{equation}\label{1.4}
|\sigma_n(v)|\leq g_n,\quad \sigma_n(v)v_n+g_n |v_n|=0,
\end{equation}
where and in what follows $v_n=v\cdot n$.

Till now, for the Stokes and Navier-Stokes problems with friction
type boundary conditions rather simple cases are studied. More
clearly, one deal with problems with the Dirichlet boundary
condition on a portion of boundary and one of friction type
conditions on other portion.

In \cite{f1} existence of solutions to the steady Stokes and
Navier-stokes equations with the homogeneous Dirichlet boundary
condition on a portion of boundary and leak  or threshold slip
boundary condition on other portion is studied. Also,
\cite{f2}-\cite{f4} concerned with the steady or non-steady  Stokes
equations with the homogeneous Dirichlet boundary condition and leak
boundary condition.

 When a portion of boundary with Dirichlet
boundary condition and other moving portion where nonlinear slip
occurs are separated, existence, uniqueness and continuous
dependence on the data are studied for the steady Stokes equations in
\cite{r} and for the steady Navier-Stokes equations in \cite{rt}. In
\cite{sa} when a portion of boundary with Dirichlet boundary
condition and another portion with slip condition are separated, existence of
strong solution to the steady Stokes equations is studied. In
\cite{sf} when a portion with homogeneous Dirichlet
boundary condition and other portion with nonlinear boundary
condition are separated, for the steady Stokes equations a relation
between a regularized problem and the original problem, regularity
of solution are studied.

In \cite{rt} for the steady Navier-Stokes equations, existence,
uniqueness and continuous dependence on the data are studied when a
portion of boundary with Dirichlet boundary condition and another
moving portion where nonlinear slip occurs are separated. In
\cite{s} local unique existence of solution to the steady
Navier-Stokes problem with homogeneous Dirichlet boundary condition
and one of friction boundary conditions is studied.  In \cite{al1}
existence and uniqueness of solution to the steady rotating
Navier-Stokes equations are studied when boundary consists of a portion with
homogeneous Dirichlet boundary condition and other portions where
there is flow and threshold slip. In \cite{ll3} under similar
boundary condition the steady Navier-Stokes problem is studied.

 In \cite{all} existence of weak solution and local existence of a strong solution to the non-steady Navier-Stokes problem are
studied when boundary consists of a portion with homogeneous
Dirichlet boundary condition and another portion with leak
condition. In \cite{k2} existence of a strong solution to the
non-steady Navier-Stokes equation is studied when boundary consists
of a portion with homogeneous Dirichlet boundary condition and
another portion with  nonlinear slip or leak condition.

For other kinds of non-steady fluid equations with friction slip boundary conditions and Dirichlet condition, refer to \cite{bbp1}, \cite{bbp2}, \cite{bl} and \cite{dr}.
   Numerical solution methods are studied  for the Stokes and Navier-Stokes problems with friction boundary conditions.
For the 2-D steady Stokes problems refer to  \cite{ags},
  \cite{k1}, \cite{ll1}, \cite{ll2},  \cite{ll4} and for the 3-D steady Stokes problems \cite{k3}. For the 2-D steady Navier-Stokes problem refer to \cite{a}, \cite{la2}, \cite{la4} and \cite{la3}.
For the 2-D non-steady Navier-Stokes problem refer to
\cite{la1}.

In practice we deal with mixture of some kinds of boundary
 conditions.  Especially, when there is flux through a portion of boundary,
 we can deal with pressure boundary conditions. There are many
 papers dealing with the total pressure
(Bernoulli's pressure) $\f{1}{2} |v|^2+p$ (cf. \cite{cmp},
\cite{cppt}) or static pressure $p$ (cf. \cite{act}, \cite{m1}). It
is also known that the total stress $\sigma^t(v,p)$ on the boundary
is a natural boundary condition, where $\sigma^t(v,p)=S^tn, $ and
total stress tensor $S^t$ is one with components
$S^t_{ij}=-(p+\f{1}{2} |v|^2)\delta_{ij}+2\nu\varepsilon_{ij}(v).$
(see \cite{fmn1}, \cite{fmn}).

Also, in practice we deal with one-sided leak of fluid. The condition
\eqref{1.4} means that according to direction of normal stress,
fluid penetrates out or into through boundary. If the fluid can only
leak out through boundary when $-\sigma_n(v)$ is same with a
threshold $g_{+n}(>0)$, then we can describe that by
\begin{equation}\label{1.5}
v_n\geq 0,\quad\sigma_n(v)+g_{+n}\geq 0,\quad
(\sigma_n(v)+g_{+n})v_n=0.
\end{equation}
In contrast, if the fluid can only leak into through boundary when
 $-\sigma_n(v)$ is same with a threshold $-g_{-n}(g_{-n}>0)$, then we can
describe that by
\begin{equation}\label{1.6}
v_n\leq 0,\quad\sigma_n(v)-g_{-n}\leq 0,\quad
(\sigma_n(v)-g_{-n})v_n=0.
\end{equation}
For one-sided flow condition depending on a threshold of total
pressure refer to \cite{c}. For similar one-sided boundary
conditions of elasticity refer to \cite{kl}, Section 5.4.1, ch. 3 in
\cite{dl}.

  In the present paper, we are
concerned with the the systems \eqref{1.1} and \eqref{1.2} with mixed
boundary conditions involving Tresca slip condition \eqref{1.3},
leak boundary condition \eqref{1.4}, one-sided leak boundary
conditions \eqref{1.5} and \eqref{1.6}, velocity, static pressure,
rotation, stress and normal derivative of velocity together. And
also without discussing whether static pressure or total pressure
(correspondingly stress or total stress) is suitable for real phenomena
which is over our knowledge, we consider the problems with total pressure and
total stress instead of static pressure and stress.
Relying on the result in \cite{kc}, we reflect all these boundary
conditions into variational formulations of problems. Overcoming
difficulty from one-sided leak boundary conditions, we get
variational inequalities equivalent to the variational formulation
for the problems. We study some variational inequalities concerned
with the Navier-Stokes problems. Using the results for the variational inequalities,
we prove existence, uniqueness and estimates of weak solutions to
the Navier-Stokes problems with such boundary conditions. Also using
the previous results for elliptic variational inequality, we get
some results for the Stokes problem with such boundary conditions.

This paper consists of 5 sections.

In Section 2, some previous results for variational formulation of
our problems are stated. Also, three problems to study are
described. For the Navier-Stokes equations, according to the
pressure or the total pressure (correspondingly stress or the total
stress) two problems are distinguished.

 In Section 3, for the stationary Navier-Stokes and Stokes problems with mixture of 11 kinds of boundary
conditions we have the variational formulations which consist of
five formulae with five unknown functions, that is, using velocity,
tangent stress on slip surface, normal stress on leak surface,
normal stresses on one-sided leak surfaces together as unknown
functions. Except friction type conditions, other boundary
conditions  are reflected in a variational equation as
usual(Problems I-VE, II-VE, III-VE). When the solution smooth
enough, these variational formulations are equivalent to the
original PDE problems(Theorems \ref{t3.0}, \ref{t3.2.0}). Then, we
get variational inequalities equivalent to the variational
formulations above, which have one unknown
function-velocity(Theorems \ref{t3.1}, \ref{t3.2}). In proof of
equivalence, to overcome difficulties from the one-sided leak
conditions Lemma \ref{l3.1} is used.

In Section 4 we study 3 kinds of variational inequality which are
for the problems in Section 3. With an exception \cite{s} studying
local unique existence, in all previous papers dealing with friction
boundary conditions one approximate the functionals in the
considering variational inequalities with smooth one resulting to
study of operator equation and it's convergence.  Owing to the
one-sided leak conditions such approximation for our problem may be
complicated. Without such approximation we first get existence,
uniqueness and estimates of solutions to the variational
inequalities(Theorems \ref{t4.1}, \ref{t4.2}). In addition, for a
special case excluding flux through boundary we also show
approximation way of the functional(Theorem \ref{t4.3}).

In Section 5, relying the results in Section 4, we study existence,
uniqueness and estimates of solutions to the Navie-Stokes problems
with 11 kinds of boundary condition. For the Navier-Stokes problem
with boundary condition \eqref{2.8}, which is including static pressure and stress,  local unique existence is
proved(Theorem \ref{t5.1}). For the Navier-Stokes problem with
boundary condition \eqref{2.9}, which is including total static pressure and total stress, existence and estimate of solutions
are proved(Theorem \ref{t5.2}).  For a special case of the
Navier-Stokes problem with boundary condition  \eqref{2.8} in which
there is no any flux across boundary except $\Gamma_1, \Gamma_8$,
existence and estimate of solutions are proved(Theorem \ref{t5.2.0},
\ref{t5.3}). Also, relying the previous results in elliptic
variational inequality, we study unique existence, an estimate and
continuous dependence on data of solutions to the Stokes problem
with the boundary condition \eqref{2.8}(Theorem \ref{t5.4}).

 Throughout this paper we will use the
following notation.

Let $\Omega$ be a connected bounded open subset of $ R^l,\,l=2,3.$
$\partial\Omega\in C^{0,1}$, $\partial\Omega=\cup_{i=1}^{11}
\overline{\Gamma}_i$,\,$\Gamma_i\cap\Gamma_j=\varnothing$ for $i\neq
j$, $\Gamma_i=\bigcup_j\Gamma_{ij}$, where $\Gamma_{ij}$ are
connected open subsets of $\partial\Omega$ and $\Gamma_{ij}\in C^2$
for $i=2,3,7$ and $\Gamma_{ij}\in C^1$ for others.
 When $X$ is a Banach space,  $\mathbf{X}=X^l$. Let $W^k_\alpha(\Omega)$ be
Sobolev spaces, $H^1(\Omega)=W^1_2(\Omega)$, and so
$\mathbf{H}^1(\Omega)=\{H^1(\Omega)\}^l$. Let $0_X$ be the zero
element of space $X$ and $\mathscr{O}_M(0_X)$ be $M$-neighborhood of
$0_X$ in space $X$. Compact continuous imbedding of a space $X$ into
a space $Y$ is denoted by $X\hookrightarrow\hookrightarrow Y$.

An inner product and norm in the space $\mathbf{L}_2(\Omega)$ are, respectively,
denoted by $(\cdot\,,\cdot)$ and $\|\cdot\|$; and
$\langle\cdot\,,\cdot\rangle$ means the duality pairing between  a
Sobolev space $X$ and its dual one. Also,
$(\cdot\,,\cdot)_{\Gamma_i}$ is an inner product in the
$\mathbf{L}_2(\Gamma_i)$ or $L_2(\Gamma_i)$; and
$\langle\cdot\,,\cdot\rangle_{\Gamma_i}$ means the duality pairing
between $\mathbf{H}^{\frac{1}{2}}(\Gamma_i)$ and
$\mathbf{H}^{-\frac{1}{2}}(\Gamma_i)$ or between
$H^{\frac{1}{2}}(\Gamma_i)$ and $H^{-\frac{1}{2}}(\Gamma_i)$. The
inner product and norms in $R^l$, respectively, are denoted by
$(\cdot\, , \cdot)_{R^l}$ and $|\cdot|$.  Sometimes the inner
product between $a$ and $b$ in $R^l$ is denoted by $a\cdot b$.
For convenience, in the case that $l=2,\, y=(y_1(x_1,x_2),\,
y_2(x_1,x_2))$ is identified with $\bar{y}=(y_1,\,y_2,\,0)$, and so
$\text{rot}\,y=\text{rot}\,\bar{y}$. Thus, for $y=(y_1,\,y_2)$ and
$v=( v_1,\,v_2),\, \text{rot}\,y\times v $ is the 2-D vector
consisted of
the first two components of $\mbox{rot}\,\bar{y}\times \bar{v}$.

Let $n(x)$ and
$\tau(x)$ be, respectively, outward normal and tangent unit vectors
 at $x$ in $\partial\Omega$. When for $u\in H^1(\Omega)$ such that $u_\tau=0$ on $\Gamma_i$, sometimes  for convenience we use notation $u|_{\Gamma_i}$ instead $u_n|_{\Gamma_i}.$
If when $f\in H^{-1/2}(\Gamma_i),$ $\left< f,
w\right>_{\Gamma_i}\geq 0\,\, (\leq 0)\,\,\forall w\in
C_0^\infty(\Gamma_i)$ with $w\geq 0 $, then we denoted by $f\geq
0\,\,( \leq 0)$.

\section{Preliminary and problems}
\setcounter{equation}{0}

 Let $\Gamma$ be a surface (curve for $l=2$) of $C^2$ and $v$ be  a vector field of $C^2$ on a domain of $R^l$ near $\Gamma$.
In this paper the surfaces concerned by us are pieces of boundary of
3-D or 2-D bounded connected domains, and so we can assume the
surfaces are oriented.
\begin{theorem}\label{t2.1} (Theorem 2.1 in \cite{kc})
Suppose that $v\cdot n|_{\Gamma}=0$. Then, on the surface $\Gamma$
the following holds.
\begin{equation}\label{2.1}
\left(\varepsilon(v) n,
\tau\right)_{R^l}=\f{1}{2}(\mbox{rot}\,v\times n,
\tau)_{R^l}-(S\tilde{v},\tilde{\tau})_{R^{l-1}},
\end{equation}
\begin{equation}\label{2.2}
(\mbox{rot}\,v\times n, \tau)_{R^l}=\left(\f{\p v}{\p n},
\tau\right)_{R^l}+(S\tilde{v},\tilde{\tau})_{R^{l-1}},
\end{equation}
\begin{equation}\label{2.3}
\left(\varepsilon(v) n, \tau\right)_{R^l}=\f{1}{2}\left(\f{\p v}{\p
n}, \tau\right)_{R^l}-\f{1}{2} (S\tilde{v},\tilde{\tau})_{R^{l-1}},
\end{equation}
where $\varepsilon(v)$ denotes the matrix with the components
$\varepsilon_{ij}(v)$, $S$ is the shape operator of the surface
$\Gamma$ (the matrix (A.1) in \cite{kc}) for $l=3$ and the curvature
of $\Gamma$ for $l=2$, and $\tilde{v}, \tilde{\tau}$ are expressions
of the vectors $v, \tau$ in a local curvilinear coordinates on
$\Gamma$.
\end{theorem}

\begin{remark}\label{r2.1}  Assuming $\Gamma$ be a surface of $C^2$, let us introduce a local curvilinear
coordinate system on $\Gamma$ which is orthogonal at all points each
other. Then, the shape operator $S$ is expressed by the following
matrix
\[
 S=
\begin{pmatrix}
L & K\\
M & N
\end{pmatrix},
\]
where
\[
\begin{aligned}
L=\left(e_1, \f{\p n}{\p e_1}\right)_{R^l}, \,\,K=\left(e_2, \f{\p
n}{\p e_1}\right)_{R^l},\,\, M=\left(e_1, \f{\p n}{\p
e_2}\right)_{R^l}, \,\, N=\left(e_2, \f{\p n}{\p e_2}\right)_{R^l}
\end{aligned}
\]
and $e_i, i=1,2,$ are unit vector of the local coordinate system.
The bilinear form $(S\tilde{v},\tilde{u})_{R^{l-1}}$ for vector
$u,v$ tangent to the surface is independent from curvilinear
coordinate system which is orthogonal at all points each other(cf.
Appendix in \cite{kc}).
\end{remark}

\begin{theorem}\label{t2.2} (Theorem 2.2 in \cite{kc})
On the surface $\Gamma$ the following holds.
\begin{equation}\label{2.4}
\left(\varepsilon(v) n, n\right)_{R^l}=\left(\f{\p v}{\p n},
n\right)_{R^l}.
\end{equation} If $v\cdot
\tau|_{\Gamma}=0$, then
\begin{equation}\label{2.5}
\left(\varepsilon(v) n, n\right)_{R^l}=\left(\f{\p v}{\p n},
n\right)_{R^l}=-(k(x)v,n)_{R^l}-
\mbox{div}_{\Gamma}\,v_{\tau}+\mbox{div}\,v,
\end{equation}
where $k(x)=\mbox {div}\,n(x)$, $v_{\tau}$ is the tangential
component of $v$ and $div_\Gamma$ is the divergence of a tangential
vector field in the tangential coordinate system on $\Gamma$.
\end{theorem}
\begin{definition}\label{d2.1} (Definition A.2 in \cite{kc})
If a piece of boundary on a neighborhood of $x\in \p \Omega$ is on
the opposite (same) side of the outward normal vector with respect
to tangent plane (line for l=2) at $x$ or coincides with the tangent
plane, then piece of the boundary called convex (concave) at $x$. If
at all $x\in \Gamma\subset \p \Omega$ the boundary convex (concave),
then $\Gamma$ called convex (concave).
\end{definition}

\begin{lemma}\label{l2.1} (Lemma A.3 in \cite{kc})
If $\Gamma_{ij}$  are convex (concave), then quadratic forms
$(S\tilde{v},\tilde{v})|_{\Gamma_i}$ and $(k(x)v,v)_{\Gamma_i}$ are
positive (negative).
\end{lemma}

\begin{definition}\label{d2.2} A functional $f:X\rightarrow \overline{R}\equiv R\cup +\infty$ is said to be proper if it is not
identically equal to $\infty$. If $f(x)\in (-\infty, +\infty)\,\,\,
\forall x\in X$, then it is said to be finite.
\end{definition}

We are concerned the problems I and II for the Navier-Stokes
equations
\begin{equation}\label{2.7}
 -\nu\Delta v+(v\cdot\nabla)v+\nabla p=f,\,\,\,\nabla \cdot v=0\quad \mbox{in}\,\,
 \Omega,
\end{equation}
which are distinguished according to boundary conditions. Problem I
is one with the boundary conditions
\begin{equation}\label{2.8}
\begin{aligned}
&(1)\quad v|_{\Gamma_1}=h_1,\\
&(2)\quad v_\tau|_{\Gamma_2}=0,\,\,-p|_{\Gamma_2}=\phi_2,\\
&(3)\quad v_n|_{\Gamma_3}=0,\,\,\mbox{rot}\,v\times n|_{\Gamma_3}=\phi_3/\nu,\\
&(4)\quad v_\tau|_{\Gamma_4}=h_4,\,\,(-p+2\nu\varepsilon_{nn}(v))|_{\Gamma_4}=\phi_4,\\
&(5)\quad v_n|_{\Gamma_5}=h_5,\,\,2(\nu\varepsilon_{n\tau}(v)+\alpha v_\tau)|_{\Gamma_5}=\phi_5,\,\,\,\alpha: \text{a matrix},\\
&(6)\quad (-pn+2\nu\varepsilon_n(v))|_{\Gamma_6}=\phi_6,\\
&(7)\quad v_\tau|_{\Gamma_7}=0,\,(-p+\nu\f{\p v}{\p n}\cdot
n)|_{\Gamma_7}=\phi_7,\\
&(8)\quad v_n|_{\Gamma_8}=h_8,\,|\sigma_\tau(v)|\leq g_\tau,\,
\sigma_\tau(v)\cdot v_\tau+g_\tau |v_\tau|=0\,\,
\mbox{on}\,\, \Gamma_8, \\
&(9)\quad v_\tau|_{\Gamma_9}=h_9,\,|\sigma_n(v)|\leq
g_n,\,\sigma_n(v) v_n+g_n |v_n|=0\,\, \mbox{on}\,\, \Gamma_9,
\\
&(10)\quad v_\tau=0,\,\, v_n\geq 0,\quad\sigma_n(v)+g_{+n}\geq
0,\quad (\sigma_n(v)+g_{+n})v_n=0
  \,\,
\mbox{on}\,\, \Gamma_{10},
\\
&(11)\quad v_\tau=0,\,v_n\leq 0,\quad\sigma_n(v)-g_{-n}\leq 0,\quad
(\sigma_n(v)-g_{-n})v_n=0\,\, \mbox{on}\,\, \Gamma_{11},
\end{aligned}
\end{equation}
and Problem II is one with the conditions
\begin{equation}\label{2.9}
\begin{aligned}
&(1)\quad v|_{\Gamma_1}=h_1,\\
&(2)\quad v_\tau|_{\Gamma_2}=0,\,\,  -(p+\f{1}{2}|v|^2)|_{\Gamma_2}=\phi_2,\\
&(3)\quad v_n|_{\Gamma_3}=0,\,\,\mbox{rot}\,v\times n|_{\Gamma_3}=\phi_3/\nu,\\
&(4)\quad v_\tau|_{\Gamma_4}=h_4,\,\,(-p-\f{1}{2}|v|^2+2\nu\varepsilon_{nn}(v))|_{\Gamma_4}=\phi_4,\\
&(5)\quad v_n|_{\Gamma_5}=h_5,\,\,2(\nu\varepsilon_{n\tau}(v)+\alpha v_\tau)|_{\Gamma_5}=\phi_5,\,\,\,\alpha: \text{a matrix},\\
&(6)\quad
(-pn-\f{1}{2}|v|^2n+2\nu\varepsilon_n(v))|_{\Gamma_6}=\phi_6,\\
&(7)\quad v_\tau|_{\Gamma_7}=0,\,(-p-\f{1}{2}|v|^2+\nu\f{\p v}{\p
n}\cdot n)|_{\Gamma_7}=\phi_7,\\
&(8)\quad v_n|_{\Gamma_8}=h_8,\,|\sigma_\tau^t(v)|\leq g_\tau,\,
\sigma_\tau^t(v)\cdot v_\tau+g_\tau |v_\tau|=0\,\,
\mbox{on}\,\, \Gamma_8, \\
&(9)\quad v_\tau|_{\Gamma_9}=h_9,\,|\sigma_n^t(v)|\leq g_n,\,
\sigma_n^t(v) v_n+g_n |v_n|=0\,\, \mbox{on}\,\, \Gamma_9,
\\
&(10)\quad v_\tau=0,\,\, v_n\geq 0,\quad\sigma^t_n(v)+g_{+n}\geq
0,\quad (\sigma^t_n(v)+g_{+n})v_n=0
  \,\,
\mbox{on}\,\, \Gamma_{10},
\\
&(11)\quad v_\tau=0,\,v_n\leq 0,\quad\sigma^t_n(v)-g_{-n}\leq
0,\quad (\sigma^t_n(v)-g_{-n})v_n=0\,\, \mbox{on}\,\, \Gamma_{11},
\end{aligned}
\end{equation}
where $\varepsilon_n(v)=\varepsilon(v)n$,
$\varepsilon_{nn}(v)=(\varepsilon(v)n,n)_{R^l}$,
$\varepsilon_{n\tau}(v)=\varepsilon(v)n-\varepsilon_{nn}(v)n$ and
$h_i, \phi_i, \alpha_{ij}$(components of matrix $\alpha)$ are given
functions or vectors of functions. And $\sigma^t_n$ is the normal
component of total stress on surface, that is,
$\sigma^t_n=\sigma^t\cdot n.$ Also,
$\sigma^t_\tau(v,p)=\sigma^t(v,p)-\sigma^t_n(v,p)n$ and
 $g_\tau\in L^2(\Gamma_8)$, $g_n\in L^2(\Gamma_9)$, $g_{+n}\in L^2(\Gamma_{10})$, $g_{-n}\in L^2(\Gamma_{11})$, $g_\tau>0$,
$g_n>0$, $g_{+n}>0$, $g_{-n}>0$, at a.e.

For Problem II the static pressure $p$ and stress in the boundary
conditions for Problem I are changed with the total pressure and the
total stress. Note
$$\sigma_\tau(v,p)=\sigma^t_\tau(v,p)=2\nu\varepsilon_{n\tau}(v).$$

We also consider the Stokes equations
\begin{equation}\label{2.6}
\begin{aligned}
 -\nu\Delta v+\nabla p=f, \,\,\,
 \nabla \cdot v=0\quad \mbox{in}\,\, \Omega
 \end{aligned}
\end{equation}
with the boundary conditions \eqref{2.8}, which is Problem III.

\section{Variational formulations and equivalent variational inequalities}
\setcounter{equation}{0}

In this section we give variational formulations for Problems I, II, III above and get variational inequalities equivalent to the formulations.

Let
$$
\begin{aligned}
&\mathbf{V}(\Omega)=\{u\in \mathbf{H}^1(\Omega):\mbox{div}\,u=0,\,
u|_{\Gamma_1}=0,\,
u_\tau|_{(\Gamma_2\cup\Gamma_4\cup\Gamma_7\cup\Gamma_9\cup\Gamma_{10}\cup\Gamma_{11})}=0,\,
u_n|_{(\Gamma_3\cup\Gamma_5\cup\Gamma_8)}=0\},\\&
\mathbf{V}_{\Gamma237}(\Omega)=\{u\in
\mathbf{H}^1(\Omega):\mbox{div}\,u=0,\,
u_\tau|_{(\Gamma_2\cup\Gamma_7)}=0,\, u_n|_{\Gamma_3}=0\},
\end{aligned}
$$
and
$$
K(\Omega)=\{u\in\mathbf{V}(\Omega): u_n|_{\Gamma_{10}}\geq
0,\,u_n|_{\Gamma_{11}}\leq 0\}.
$$

By Theorem \ref{2.1} and \ref{2.2} for $v\in
\mathbf{H}^2(\Omega)\cap\mathbf{V}_{\Gamma237}(\Omega)$ and $u\in
\mathbf{V}(\Omega)$
\begin{equation}\label{3.1}
\begin{aligned}
 -(\Delta v,
u)&=2(\varepsilon(v), \varepsilon(u))-2(\varepsilon(v)n,
u)_{\cup_{i=2}^{11}\Gamma_i}\\& =2(\varepsilon(v),
\varepsilon(u))+2(k(x)v,u)_{\Gamma_2}-(\mbox{rot}\,v\times n,
u)_{\Gamma_3}+2(S\tilde{v},\tilde{u})_{\Gamma_3}\\&\quad-2(\varepsilon_{nn}(v),u_n)_{\Gamma_4}-2(\varepsilon_{n\tau}(v),u)_{\Gamma_5}-2(\varepsilon_n(v),u)_{\Gamma_6}-\left(\f{\p
v}{\p
n},u\right)_{\Gamma_7}+(k(x)v,u)_{\Gamma_7}\\&\quad-2(\varepsilon_{n\tau}(v),u)_{\Gamma_8}-2(\varepsilon_{nn}(v),u_n)_{\Gamma_9}-2(\varepsilon_{nn}(v),u_n)_{\Gamma_{10}}-2(\varepsilon_{nn}(v),u_n)_{\Gamma_{11}}.
\end{aligned}
\end{equation}
Also, for $p\in H^1(\Omega)$ and $u\in \mathbf{V}(\Omega)$ we have
\begin{equation}\label{3.2}
\begin{aligned}
(\nabla p, u)&=(p, u_n)_{\cup_{i=2}^{11}\Gamma_i}=(p, u_n)_{\Gamma_2}+(p, u_n)_{\Gamma_4\cup\Gamma_7\cup\Gamma_9\cup\Gamma_{10}\cup\Gamma_{11}}+(pn,
u)_{\Gamma_6},
\end{aligned}
\end{equation}
where $u_n\mid_{\Gamma_3\cup\Gamma_5\cup\Gamma_8}=0$ was used.\\

We assume that the following holds.
\begin{assumption}\label{a3.1}
1) There exists a function $U\in\mathbf{H}^1(\Omega)$ such that
$$
\begin{aligned}
&\mbox{div}\,U=0,U|_{\Gamma_1}=h_1,
U_\tau|_{(\Gamma_2\cup\Gamma_7)}=0,
U_n|_{\Gamma_3}=0, U_\tau|_{\Gamma_4}=h_4,\\& U_n|_{\Gamma_5}=h_5,
U|_{\Gamma_8}=h_8
n,\,U|_{\Gamma_9}=h_9,\,U|_{\Gamma_{10}}=0,\,U|_{\Gamma_{11}}=0.
\end{aligned}
$$

2) $f\in\mathbf{V}(\Omega)^*,\,\, \phi_i\in
H^{-\f{1}{2}}(\Gamma_i),i=2,4,7,\,\,
\phi_i\in\mathbf{H}^{-\f{1}{2}}(\Gamma_i),i=3,5,6,\,\,\alpha_{ij}\in
L_\infty(\Gamma_5)$, and $\Gamma_1\neq\varnothing$.

3) If $\Gamma_i,$ where $i$ is 10 or 11, is nonempty, then at least
one of $\{\Gamma_j:j\in\{2,4,7,9-11\}\backslash i\}$ is nonempty and
there exist a diffeomorphisms in $C^1$ between $\Gamma_i$ and
$\Gamma_j$.
\end{assumption}

Having in mind Assumption \ref{a3.1} and putting $v=w+U$, by
\eqref{3.1}, \eqref{3.2} we can see that smooth solutions $v$ of
problem \eqref{2.7}, \eqref{2.8} satisfy the following.

\begin{equation}\label{3.3}
\left\{
\begin{aligned}
&v-U=w\in K(\Omega),\\&
 2\nu(\varepsilon(w),\varepsilon(u))+\langle (w\cdot\nabla)
w,u\rangle+\langle(U\cdot
\nabla)w,u\rangle+\langle(w\cdot\nabla)U,u\rangle\\&\hspace*{1cm}+2\nu(k(x)w,u)_{\Gamma_2}+2\nu(S\tilde{w},\tilde{u})_{\Gamma_3}
+2(\alpha(x)w,u)_{\Gamma_5}+\nu(k(x)w,u)_{\Gamma_7}\\&
\hspace*{1cm}-2(\varepsilon_{n\tau}(w+U),u)_{\Gamma_8}
+(p-2\varepsilon_{nn}(w+U),u_n)_{\Gamma_9\cup\Gamma_{10}\cup\Gamma_{11}}\\&
\hspace*{.5cm}=-2\nu(\varepsilon(U),\varepsilon(u))-\langle(U\cdot\nabla)U,u\rangle
-2\nu(k(x)U,u)_{\Gamma_2}-2\nu(S\tilde{U},\tilde{u})_{\Gamma_3}-2(\alpha(x)U,u)_{\Gamma_5}\\&
\hspace*{1cm}-\nu(k(x)U,u)_{\Gamma_7} +\langle
f,u\rangle+\sum_{i=2,4,7}\langle\phi_i,u_n\rangle_{\Gamma_i}+\sum_{i=3,5,6}\langle\phi_i,u\rangle_{\Gamma_i}
\quad \forall u\in \textbf{V}(\Omega),\\
&|\sigma_\tau(v)|\leq g_\tau,\, \sigma_\tau(v)\cdot v_\tau+g_\tau
|v_\tau|=0\quad \mbox{on}\,\, \Gamma_8,\\
&|\sigma_n(v)|\leq g_n,\, \sigma_n(v) v_n+g_n |v_n|=0\quad
\mbox{on}\,\, \Gamma_9,\\
 &
 \sigma_n(v)+g_{+n}\geq 0,\quad
(\sigma_n(v)+g_{+n}) v_n=0\quad \mbox{on}\,\, \Gamma_{10},
\\
&\sigma_n(v)-g_{-n}\leq 0,\quad (\sigma_n(v)-g_{-n}) v_n=0\quad
\mbox{on}\,\, \Gamma_{11}.
\end{aligned}
\right.
\end{equation}
Define $a_{01}(\cdot, \cdot), a_{11}(\cdot,\cdot,\cdot)$ and $F_1\in
V^*$ by
\begin{equation}\label{3.4}
\begin{aligned}
& a_{01}(w,u)=2\nu(\varepsilon(w),\varepsilon(u))+\langle(U\cdot
\nabla)w,u\rangle+\langle(w\cdot\nabla)U,u\rangle+2\nu(k(x)w,u)_{\Gamma_2}\\&
\hspace{1.6cm}+2\nu(S\tilde{w},\tilde{u})_{\Gamma_3}+2(\alpha(x)w,u)_{\Gamma_5}+\nu(k(x)w,u)_{\Gamma_7}\quad
\forall w,u\in \mathbf{V}(\Omega),
\\&a_{11}(w,u,v)=\langle (w\cdot\nabla)
u,v\rangle\quad \forall w,u,v\in \mathbf{V}(\Omega),\\
& \langle
F_1,u\rangle=-2\nu(\varepsilon(U),\varepsilon(u))-\langle(U\cdot\nabla)U,u\rangle-2\nu(k(x)U,u)_{\Gamma_2}-2\nu(S\tilde{U},\tilde{u})_{\Gamma_3}\\&
\hspace*{1.3cm}-2(\alpha(x)U,u)_{\Gamma_5}-\nu(k(x)U,u)_{\Gamma_7}
+\langle
f,u\rangle+\sum_{i=2,4,7}\langle\phi_i,u_n\rangle_{\Gamma_i}\\&\hspace*{1.3cm}+\sum_{i=3,5,6}\langle\phi_i,u\rangle_{\Gamma_i}\quad
\forall u\in \mathbf{V}(\Omega).
\end{aligned}
\end{equation}

Then, taking into account
\[
\sigma_\tau(v)=2\nu\varepsilon_{n\tau}(v),\quad
\sigma_n(v)=-p+2\nu\varepsilon_{nn}(v)
\]
and \eqref{3.3}, we introduce the following variational formulation for
problem \eqref{2.7}, \eqref{2.8}.\vspace{.1cm}

\textbf{Problem I-VE.} Find $(v, \sigma_\tau, \sigma_n, \sigma_{+n},
\sigma_{-n})\in \big(U+K(\Omega)\big)\times
\textbf{L}^2_\tau(\Gamma_8)\times L^2(\Gamma_9)\times
H^{-1/2}(\Gamma_{10})\times H^{-1/2}(\Gamma_{11})$ such that
\begin{equation}\label{3.5}
\left\{\begin{aligned} &v-U=w\in K(\Omega),\\&
a_{01}(w,u)+a_{11}(w,w,u)-(\sigma_\tau,u_\tau)_{\Gamma_8}-(\sigma_n,u_n)_{\Gamma_9}\\&
\hspace{1.4cm}-\left<\sigma_{+n},u_n\right>_{\Gamma_{10}}-\left<\sigma_{-n},u_n\right>_{\Gamma_{11}}=\langle
F_1,u\rangle\quad  \forall u\in \mathbf{V}(\Omega),
\\
&|\sigma_\tau|\leq g_\tau,\, \sigma_\tau\cdot v_\tau+g_\tau
|v_\tau|=0\quad \mbox{on}\,\, \Gamma_8,\\
&|\sigma_n|\leq g_n,\, \sigma_n v_n+g_n |v_n|=0\quad \mbox{on}\,\,
\Gamma_9,
\\&\sigma_{+n}+g_{+n}\geq 0,\quad
\left<\sigma_{+n}+g_{+n}, v_n\right>_{\Gamma_{10}}=0\quad \mbox{on}\,\, \Gamma_{10},
\\
&\sigma_{-n}-g_{-n}\leq 0,\quad \left<\sigma_{-n}-g_{-n}, v_n\right>_{\Gamma_{11}}=0\quad
\mbox{on}\,\, \Gamma_{11},
\end{aligned}
\right.
\end{equation}
where $\textbf{L}^2_\tau(\Gamma_8)$ is the subspace of
$\textbf{L}^2(\Gamma_8)$  consisted of functions such that $ (u,
n)_{\textbf{L}^2(\Gamma_8)}=0.$

\begin{remark}\label{r3.1}
If $u\in H^1(\Omega),$ then $v|_{\Gamma_i}\in
H^{\f{1}{2}}(\Gamma_i),$ however if $u|_{\p\Omega}=0$ on
$O(\Gamma_i)\backslash\overline{\Gamma}_i,$ where $O(\Gamma_i)$ is
an open subset of $\p\Omega$ such that $\overline{\Gamma}_i\subset
O(\Gamma_i),$ then $u|_{\Gamma_i}\in H^{\f{1}{2}}_{00}(\Gamma_i)$
(cf. (c) of Theorem 1.5.2.3 in \cite{g1}). Since
$H^{\f{1}{2}}_{00}(\Gamma_i)\hookrightarrow
H_0^{\f{1}{2}}(\Gamma_i)$ and
$H^{\f{1}{2}}(\Gamma_i)=H_0^{\f{1}{2}}(\Gamma_i),$ (cf. Theorems
11.7 and 11.1 of ch. 1 in \cite{lm})
$$
H_{00}^{\f{1}{2}}(\Gamma_i)\hookrightarrow
H^{\f{1}{2}}(\Gamma_i)\hookrightarrow
H^{-\f{1}{2}}(\Gamma_i)\hookrightarrow
(H_{00}^{\f{1}{2}}(\Gamma_i))'.
$$
Thus, under condition $u|_{\p\Omega}=0$ on
$O(\Gamma_i)\backslash\overline{\Gamma}_i,$ for
$\phi_i\in(H^{\f{1}{2}}_{00}(\Gamma_i))'$ a dual product
$\langle\phi_i,u\rangle_{\Gamma_i}$ has meaning. But, without
knowing that $u|_{\p\Omega}=0$ on
$O(\Gamma_i)\backslash\overline{\Gamma}_i,$  for $\phi_i\in
H^{-\f{1}{2}}(\Gamma_i)$ the dual product
$\langle\phi_i,u\rangle_{\Gamma_i}$ has meaning. Therefore, under 2)
of Assumption \ref{a3.1} the dual products on $\Gamma_i$ in
\eqref{3.3} have meaning.
\end{remark}

\begin{theorem}\label{t3.0} Assume 1), 2) of Assumption \ref{a3.1}. If a
solution smooth enough ($v\in\textbf{H}^2(\Omega),f\in
\textbf{L}^2(\Omega)$), then Problem I-VE is equivalent to problem
\eqref{2.7}, \eqref{2.8}. In addition, if among
$\Gamma_i,i=2,4,6,7,9-11,$ at least one is nonempty, then $p$ of
problem \eqref{2.7}, \eqref{2.8} is unique.
\end{theorem}

$\mathbf{Proof.}$\, It is enough to prove conversion from Problem
I-E to problem \eqref{2.7}, \eqref{2.8}.

 Let $v$ is a solution smooth enough to Problem I-VE.
From  \eqref{3.4},  \eqref{3.5} we have
\begin{equation}\label{3.5.1}
\begin{aligned}
 &2\nu(\varepsilon(v),\varepsilon(u))+\langle (v\cdot\nabla)
v,u\rangle+2\nu(k(x)v,u)_{\Gamma_2}+2\nu(S\tilde{v},\tilde{u})_{\Gamma_3}
\\&
\hspace*{.3cm}+2(\alpha(x)v,u)_{\Gamma_5}+\nu(k(x)v,u)_{\Gamma_7}-(\sigma_\tau,u_\tau)_{\Gamma_8}-(\sigma_n,u_n)_{\Gamma_9}
\\&
\hspace*{.3cm}-\left<\sigma_{+n},u_n\right>_{\Gamma_{10}}-\left<\sigma_{-n},u_n\right>_{\Gamma_{11}}-\sum_{i=2,4,7}\langle\phi_i,u_n\rangle_{\Gamma_i}-\sum_{i=3,5,6}\langle\phi_i,u\rangle_{\Gamma_i}\\
&=\langle f,u\rangle \quad \forall u\in \textbf{V}(\Omega).
\end{aligned}
\end{equation}
From \eqref{3.1} we get
\begin{equation}\label{3.5.2}
\begin{aligned}
2\nu(\varepsilon(v), \varepsilon(u))=& -\nu(\Delta v,
u)-2\nu(k(x)v,u)_{\Gamma_2}+\nu(\mbox{rot}\,v\times n,
u)_{\Gamma_3}-2\nu(S\tilde{v},\tilde{u})_{\Gamma_3}\\&+2\nu(\varepsilon_{nn}(v),u\cdot
n)_{\Gamma_4}+2\nu(\varepsilon_{n\tau}(v),u)_{\Gamma_5}+2\nu(\varepsilon_n(v),u)_{\Gamma_6}+\nu\left(\f{\p
v}{\p
n},u\right)_{\Gamma_7}\\&-\nu(k(x)v,u)_{\Gamma_7}+2\nu(\varepsilon_{n\tau}(v),u)_{\Gamma_8}+2\nu(\varepsilon_{nn}(v),u)_{\Gamma_9}\\&
+2\nu(\varepsilon_{nn}(v),u)_{\Gamma_{10}}+2\nu(\varepsilon_{nn}(v),u)_{\Gamma_{11}}.
\end{aligned}
\end{equation}
From \eqref{3.5.1}, \eqref{3.5.2} we have
\begin{equation}\label{3.5.3}
\begin{aligned}
&(-\nu\Delta v+(v\cdot\nabla) v-f,u)+\nu(\mbox{rot}\,v\times n,
u)_{\Gamma_3}+2\nu(\varepsilon_{nn}(v),u\cdot
n)_{\Gamma_4}\\&+2\nu(\varepsilon_{n\tau}(v),u)_{\Gamma_5}+2(\alpha(x)v,u)_{\Gamma_5}+2\nu(\varepsilon_n(v),u)_{\Gamma_6}+\nu\left(\f{\p
v}{\p n},u\right)_{\Gamma_7}\\&
+2\nu(\varepsilon_{n\tau}(v),u)_{\Gamma_8}+2\nu(\varepsilon_{nn}(v),u_n)_{\Gamma_9}+2\nu(\varepsilon_{nn}(v),u_n)_{\Gamma_{10}}+2\nu(\varepsilon_{nn}(v),u_n)_{\Gamma_{11}}\\&
-(\sigma_\tau,u_\tau)_{\Gamma_8}-(\sigma_n,u_n)_{\Gamma_9}
-\left<\sigma_{+n},u_n\right>_{\Gamma_{10}}-\left<\sigma_{-n},u_n\right>_{\Gamma_{11}}\\&-\sum_{i=2,4,7}\langle\phi_i,u_n\rangle_{\Gamma_i}-\sum_{i=3,5,6}\langle\phi_i,u\rangle_{\Gamma_i}=0\\
\end{aligned}
\end{equation}
Taking any $u\in C_0^\infty$ with $div\, u=0$, we have
$$
(-\nu\Delta v+(v\cdot\nabla) v-f,u)=0,
$$
which implies existence of a unique $P\in H^1(\Omega)$ such that
$\int_\Omega P\,dx=0$ and
\begin{equation}\label{3.5.4}
-\nu\Delta v+(v\cdot\nabla) v-f=-\nabla P.
\end{equation}
(cf. Proposition 1.1, ch. 1 of \cite{t}).

Substituting \eqref{3.5.4} into \eqref{3.5.3}, integrating by parts
and taking into account \eqref{3.2}, we have
\begin{equation}\label{3.5.5}
\begin{aligned}
&(-P-\phi_2,u_n)_{\Gamma_2}+\nu(\mbox{rot}\,v\times n-\phi_3/\nu,
u)_{\Gamma_3}+(-P+2\nu\varepsilon_{nn}(v)-\phi_4,u_n)_{\Gamma_4}\\&+(2\nu\varepsilon_{n\tau}(v)
+\alpha(x)v_\tau-\phi_5,u)_{\Gamma_5}+(-Pn+2\nu\varepsilon_n(v)-\phi_6,u)_{\Gamma_6}\\&
+\big(-P+\nu\f{\p v}{\p n}\cdot n-\phi_7,u_n\big)_{\Gamma_7}
+(2\nu\varepsilon_{n\tau}(v)-\sigma_\tau,u)_{\Gamma_8}+(-P+2\nu\varepsilon_{nn}(v)-\sigma_n,u_n)_{\Gamma_9}
\\&+(-P+2\nu\varepsilon_{nn}(v)-\sigma_{+n},u_n)_{\Gamma_{10}}+(-P+2\nu\varepsilon_{nn}(v)-\sigma_{-n},u_n)_{\Gamma_{11}}=0,
\end{aligned}
\end{equation}
where $(v,u)_{\Gamma_5}=(v_\tau,u)_{\Gamma_5}$ and $(\nu\f{\p v}{\p
n},u\big)_{\Gamma_7}=(\nu\f{\p v}{\p n}\cdot n,u_n\big)_{\Gamma_7}$
were used.

Taking any  $u\in \textbf{V}$ such that $u_n|_{\p \Omega}=0,$
$u|_{\p \Omega}=0$ on $\p\Omega\setminus\Gamma_i$, respectively, for
$i=3,5,8,$ from \eqref{3.5.5} we get
\begin{equation}\label{3.5.6}
\begin{aligned}
&\mbox{rot}\,v\times n=\phi_3/\nu\quad\mbox{on}\,\,\Gamma_3,\\&
2\nu\varepsilon_{n\tau}(v)
+\alpha(x)v_\tau-\phi_5=0\quad\mbox{on}\,\,\Gamma_5,\\&
2\nu\varepsilon_{n\tau}(v)-\sigma_\tau=0\quad\mbox{on}\,\,\Gamma_8.
\end{aligned}
\end{equation}

If  for all $i=2,4,6,7,9-11,$ $\Gamma_i=\emptyset$, then putting
$p=P+c$, where $c$ is any constant, we get a solution $(v,p)$ to
problem \eqref{2.7}, \eqref{2.8}.

Assume that among $\Gamma_i,i=2,4,6,7,9-11,$ at least one is
nonempty. Taking any $u\in \textbf{V}$ such that $u_\tau|_{\p
\Omega}=0, u|_{\p \Omega}=0$ on $\p\Omega\setminus\Gamma_i$,
respectively, for $i=2,4,7,9-11,$ from \eqref{3.5.5} we have that
for some constants $c_i$, respectively,
\begin{equation}\label{3.5.7}
\begin{aligned}
&-P-\phi_2=c_2\quad\mbox{on}\,\,\Gamma_2,\\
&-P+2\nu\varepsilon_{nn}(v)-\phi_4=c_4\quad\mbox{on}\,\,\Gamma_4,\\&
 -P+\nu\f{\p v}{\p n}\cdot n-\phi_7=c_7\quad\mbox{on}\,\,\Gamma_7,\\
&-P+2\nu\varepsilon_{nn}(v)-\sigma_n=c_9\quad\mbox{on}\,\,\Gamma_9,
\\&-P+2\nu\varepsilon_{nn}(v)-\sigma_{+n}=c_{10}\quad\mbox{on}\,\,\Gamma_{10},\\
&-P+2\nu\varepsilon_{nn}(v)-\sigma_{-n}=c_{11}\quad\mbox{on}\,\,\Gamma_{11}.
\end{aligned}
\end{equation}
Taking any $u\in \textbf{V}$ such that $u|_{\p \Omega}=0$ on
$\p\Omega\setminus\Gamma_6$, from \eqref{3.5.5} we have that for a
constant $c_6$
$$
-Pn+2\nu\varepsilon_n(v)-\phi_6=c_6n\quad\mbox{on}\,\,\Gamma_6.
$$
Let us prove that all $c_i$ are equal to one constant $c$. For
example, assume that $\Gamma_2$ and $\Gamma_4$ are nonempty. Taking
any $u\in \textbf{V}$ such that $u|_{\p \Omega}=0$ on
$\p\Omega\setminus(\Gamma_2\cup\Gamma_4)$, from \eqref{3.5.5} we get
$$c_2\int_{\Gamma_2}u_n\,dx+c_4\int_{\Gamma_4}u_n\,dx=0,
$$
which implies $c_2=c_4=c$ since
$\int_{\Gamma_2}u_n\,dx=-\int_{\Gamma_4}u_n\,dx.$ Thus, from
\eqref{3.5.4}, \eqref{3.5.7}, we know that uniquely determined
$p=P+c$ satisfies
\begin{equation}\label{3.5.8}
-\nu\Delta v+(v\cdot\nabla)+\nabla p=f,
\end{equation}
and
\begin{equation}\label{3.5.9}
\begin{aligned}
&-p=\phi_2\quad\mbox{on}\,\,\Gamma_2,\\
&-p+2\nu\varepsilon_{nn}(v)=\phi_4\quad\mbox{on}\,\,\Gamma_4,\\&
-pn+2\nu\varepsilon_n(v)=\phi_6\quad\mbox{on}\,\,\Gamma_6,\\&
 -p+\nu\f{\p v}{\p n}\cdot n=\phi_7\quad\mbox{on}\,\,\Gamma_7,\\
&-p+2\nu\varepsilon_{nn}(v)=\sigma_n\quad\mbox{on}\,\,\Gamma_9,
\\&-p+2\nu\varepsilon_{nn}(v)=\sigma_{+n}\quad\mbox{on}\,\,\Gamma_{10},\\
&-p+2\nu\varepsilon_{nn}(v)=\sigma_{-n}\quad\mbox{on}\,\,\Gamma_{11}
\end{aligned}
\end{equation}
together.  By virtue of \eqref{3.5}, \eqref{3.5.6}, \eqref{3.5.9},
all conditions in \eqref{2.8} are satisfied. Therefore, $(v,p)$ is a
solution to problem \eqref{2.7}, \eqref{2.8}. $\square$\\

We will find a variational inequality equivalent to Problem I-VE.

Let $(v, \sigma_\tau, \sigma_n, \sigma_{+n},
\sigma_{-n})$ be a solution of Problem I-VE. From the second formula of
\eqref{3.5} subtracting the formula putted $u=w$ in the second
formula of \eqref{3.5}, we get
\begin{equation}\label{3.6}
\begin{aligned}
a_{01}(w,u-w)&+a_{11}(w,w,u-w)-(\sigma_\tau,u_\tau-w_\tau)_{\Gamma_8}-(\sigma_n,u_n-w_n)_{\Gamma_9}\\
&-\left<\sigma_{+n},u_n-w_n\right>_{\Gamma_{10}}-\left<\sigma_{-n},u_n-w_n\right>_{\Gamma_{11}}=\langle
F_1,u-w\rangle\quad  \forall u\in \mathbf{V}(\Omega).
\end{aligned}
\end{equation}

Define the functionals $j_\tau, j_n, j_+, j_-,$ respectively, by
\begin{equation}\label{3.7}
\begin{aligned}
&j_\tau(\eta)=\int_{\Gamma_8} g_\tau|\eta|\,dx\quad \forall \eta \in
\textbf{L}^2_\tau(\Gamma_8),\\
&j_n(\eta)=\int_{\Gamma_9} g_n|\eta|\,dx\quad \forall \eta \in
L^2(\Gamma_9),\\&
 j_+(\eta)=\int_{\Gamma_{10}} g_{+n}\eta\,dx\quad
\forall \eta \in L^2(\Gamma_{10}),
\\&
 j_-(\eta)=-\int_{\Gamma_{11}}
g_{-n}\eta\,dx\quad \forall \eta \in L^2(\Gamma_{11}).
\end{aligned}
\end{equation}
Since if $u\in K(\Omega)$, then
$u|_{\Gamma_8}\in\textbf{L}^2_\tau(\Gamma_8),$ $u_n|_{\Gamma_9}\in
L^2(\Gamma_9)$, $u_n|_{\Gamma_{10}}\in L^2(\Gamma_{10})$,
$u_n|_{\Gamma_{11}}\in L^2(\Gamma_{11})$, in what follows for
convenience we use the notation
\[
j_\tau(u)=j_\tau(u|_{\Gamma_8}),\,\,
 j_n(u)=j_n(u_n|_{\Gamma_9}),\,\,
 j_+(u)=j_+(u_n|_{\Gamma_{10}}),\,\,
 j_-(u)=j_-(u_n|_{\Gamma_{11}})\quad \forall u\in
K(\Omega).
\]

Define a functional $J(v)\in (\textbf{V}(\Omega)\rightarrow
\overline{R})$ by
\begin{equation}\label{3.8}
J(u)=\begin{cases}\begin{aligned}& j_\tau(u)+
 j_n(u)+
 j_+(u)+
 j_-(u)&\forall u\in
K(\Omega),\\& +\infty & \forall u\notin K(\Omega).
\end{aligned}
\end{cases}
\end{equation}
Then, $J$ is proper convex lower semi-continuous.

By Assumption \ref{a3.1}, $w_\tau=v_\tau$ on $\Gamma_8$ and
$w_n=v_n$ on $\Gamma_9\sim\Gamma_{11}$. Taking into account the fact
that $g_\tau|v_\tau|+\sigma_\tau\cdot v_\tau=0$, $|\sigma_\tau|\leq
g_\tau$, we have that
\begin{equation}\label{3.9}
\begin{aligned}
&j_\tau(u)-j_\tau(w)+\left(\sigma_\tau,
u_\tau\right)_{\Gamma_8}-\left(\sigma_\tau, w_\tau\right)_{\Gamma_8}
\\&\hspace{.5cm}=\int_{\Gamma_8}
\left(g_\tau|u_\tau|+\sigma_\tau \cdot
u_\tau\right)\,dx-\int_{\Gamma_8}
\left(g_\tau|w_\tau|+\sigma_\tau\cdot w_\tau\right)\,dx
\\&\hspace{.5cm}=\int_{\Gamma_8}
\left(g_\tau|u_\tau|+\sigma_\tau \cdot
u_\tau\right)\,dx-\int_{\Gamma_8}
\left(g_\tau|v_\tau|+\sigma_\tau\cdot v_\tau\right)\,dx\geq
0\quad\forall u\in K(\Omega).
\end{aligned}
\end{equation}
Taking into account the fact
that $g_n|v_n|+\sigma_n\cdot v_n=0$ and $|\sigma_n|\leq
g_n$, in the same way we have
\begin{equation}\label{3.10}
\begin{aligned}
 &j_n(u)-j_n(w)+\left(\sigma_n,
u_n\right)_{\Gamma_9}-\left(\sigma_n, w_n\right)_{\Gamma_9}\geq 0.
\end{aligned}
\end{equation}
Also,
\begin{equation}\label{3.11}
\begin{aligned}
 &j_+(u)-j_+(w)+\left<\sigma_{+n},
u_n\right>_{\Gamma_{10}}-\left<\sigma_{+n},
w_n\right>_{\Gamma_{10}}\\&=\left<g_{+n}+\sigma_{+n}, u_n\right>_{\Gamma_{10}}-\left<g_{+n}+\sigma_{+n}, w_n\right>_{\Gamma_{10}}\geq 0,\\
\end{aligned}
\end{equation}
where the facts that $u_n\geq 0,\,\sigma_{+n}+g_{+n}\geq 0$ and
$\left<\sigma_{+n}+g_{+n}, v_n\right>_{\Gamma_{10}}=0,$ $w_n=v_n$ on $\Gamma_{10}$ were used.
In the same way, we have
\begin{equation}\label{3.12}
\begin{aligned}
 &j_-(u)-j_-(w)+\left<\sigma_{-n},
u_n\right>_{\Gamma_{11}}-\left<\sigma_{-n},
w_n\right>_{\Gamma_{11}}\geq 0.
\end{aligned}
\end{equation}
By virtue of \eqref{3.8}-\eqref{3.12}, we have
\begin{equation}\label{3.13}
\begin{aligned}
J(u)-J(w)&\geq
-(\sigma_\tau,u_\tau-w_\tau)_{\Gamma_8}-(\sigma_n,u_n-w_n)_{\Gamma_9}\\
&-\left<\sigma_{+n},u_n-w_n\right>_{\Gamma_{10}}-\left<\sigma_{-n},u_n-w_n\right>_{\Gamma_{11}}\quad
\forall u\in \mathbf{V}.
\end{aligned}
\end{equation}
 Therefore, from \eqref{3.6} and \eqref{3.13} we get
\begin{equation}\label{3.15}
\begin{aligned}
a_{01}(w,u-w)+a_{11}(w,w,u-w)&+J(u)-J(w)\geq\langle
F_1,u-w\rangle\quad  \forall u\in \mathbf{V}(\Omega).
\end{aligned}
\end{equation}
Thus, we come to the following formulation associated with Problem I by
a variational inequality.\vspace*{.3cm}

\textbf{Problem I-VI} Find $v=w+U$ such that
\begin{equation}\label{3.16}
\begin{aligned}
&a_{01}(w,u-w)+a_{11}(w,w,u-w)+J(u)-J(w)\geq\langle
F_1,u-w\rangle\quad  \forall u\in \mathbf{V}(\Omega),
\end{aligned}
\end{equation}
where $a_{01}, a_{11}, F_1$ are in \eqref{3.4}, $U$ is in Assumption
\ref{a3.1} and $J$ is in \eqref{3.8}.\vspace*{.4cm}

To prove equivalence of Problem I-VI and Problem I-VE we need
\begin{lemma}\label{l3.1} For $\psi \in C^\infty_0(\Gamma_i), i=10, 11,$
there exists a function $\overline{u}\in \textbf{V}$ such that
\[
\overline{u}_n|_{\Gamma_i}=\psi,\quad
\|\overline{u}\|_{\textbf{V}}\leq C_i\|\psi\|_{H^{1/2}(\Gamma_i)},
\]
where $C_i$ are independent of $\psi$.
\end{lemma}

$\mathbf{Proof.}$\, By 3) of Assumption \ref{a3.1} if
$\Gamma_{10}\cup\Gamma_{11}\neq\varnothing$, then, for example,
$\Gamma_2\neq\varnothing$ and there exists a diffeomorphysm
$y=f_i(x)\in C^1$   from $\Gamma_i$ onto $\Gamma_2.$ Define
$\varphi(y)$ at point $y\in \Gamma_2$ corresponding to point $x\in
\Gamma_i$ by $\varphi(y)=\f{1}{Df(x)}\psi(f_i^{-1}(y)),$
 where $Df_i(x)$ is Jacobian of the transformation $f_i$. Then,
\begin{equation}\label{3.16.1}
\int_{\Gamma_2} \varphi(y)\,dy=\int_{\Gamma_i}\f{1}{Df_i(x)}\psi(f_i^{-1}(y))Df_i(x)\,dx=\int_{\Gamma_i}\psi(x)\,dx,
\end{equation}
and
\begin{equation}\label{3.16.2}
\begin{aligned}
\|\varphi(y)\|_{H^{\f{1}{2}}(\Gamma_2)}&\leq \left\|\f{1}{Df_i(x)}\right\|_{C(\overline{\Gamma_i})}\|\psi(x)\|_{H^{\f{1}{2}}(\Gamma_i)}\leq
c_i\|\psi(x)\|_{H^{\f{1}{2}}(\Gamma_i)}.
\end{aligned}
\end{equation}

When $\psi \in C^\infty_0(\Gamma_i),$  define a function
$\overline{\phi}\in \textbf{H}^{1/2}(\p\Omega)$ on $\p\Omega$ as
follows.
\[
\begin{aligned}
&  \overline{\phi}\times n|_{\Gamma_2\cup
\Gamma_i}=0,\,\,\overline{\phi}_n|_{\Gamma_2}=-\varphi,\,\,\overline{\phi}_n|_{\Gamma_{10}}=\psi,\,\,
\overline{\phi}|_{(\cup_{i=1, 3-9,11}\Gamma_i)}=0.
\end{aligned}
\]
Thus, by \eqref{3.16.1} $\int_{\p\Omega} \overline{\phi}_n\,ds=0.$
Then, there exists a solution $\overline{u}\in
\textbf{W}^{1,2}(\Omega)$ to the Stokes problem
 \[
\left\{\begin{aligned} & -\Delta u+\nabla p=0, \\& \mbox{div}\, u=0,
\\& u|_{\p\Omega}=\overline{\phi}
\end{aligned}
\right.
\]
and
\[
\|\overline{u}\|_{\textbf{V}(\Omega)}\leq c\|\overline{\phi}\|_{\textbf{H}^{1/2}(\p\Omega)}.
\]
(cf. Theorem IV.1.1 in \cite{g}). Taking into account \eqref{3.16.2}, we come to the asserted estimation with $C_i=1+c_i$. Thus $\overline{u}$ is the
asserted function. $\square$ \vspace{.2cm}

 Problem I-VE and Problem I-VI are equivalent in
the following sense.

\begin{theorem}\label{t3.1} If $(v, \sigma_\tau, \sigma_n, \sigma_{+n},
\sigma_{-n})$ is a solution to Problem I-VE, then $v$ is a solution
to Problem I-VI. Inversely, if $v$ is a solution to Problem I-VI,
then there exist $\sigma_\tau, \sigma_n, \sigma_{+n}, \sigma_{-n}$
such that $(v, \sigma_\tau, \sigma_n, \sigma_{+n}, \sigma_{-n})$  is
a solution to Problem I-VE.
\end{theorem}

$\mathbf{Proof.}$\, We already showed that if $(v, \sigma_\tau,
\sigma_n, \sigma_{+n}, \sigma_{-n})$ is a solution to Problem I-VE,
then $v$ is a solution to Problem I-VI. Thus, it is enough to prove
that if $v$ is a solution to Problem I-VI, then there exist
$\sigma_\tau, \sigma_n, \sigma_{+n}, \sigma_{-n}$ such that $(v,
\sigma_\tau, \sigma_n, \sigma_{+n}, \sigma_{-n})$  is a solution to
Problem I-VE.

Since the functional $J$ is proper, from \eqref{3.16} we have
\begin{equation}\label{3.17}
v-U=w\in K(\Omega)
\end{equation}
because if $w\notin K(\Omega),$ then the left hand side of \eqref{3.16} is $-\infty$ which is a contradiction to the fact that the right hand side is finite.

Let $\psi\in \mathbf{V}_{8-11}(\Omega)\equiv\{u\in
\mathbf{V}(\Omega) :
u|_{\Gamma_8\cup\Gamma_9\cup\Gamma_{10}\cup\Gamma_{11}}=0\}\left(\subset
K(\Omega)\right)$. Putting $u=w+\psi, u=w-\psi$ and taking into
account
\[
j_\tau(w)=j_\tau(w+\psi),\,j_n(w)=j_n(w+\psi),\,j_+(w)=j_+(w+\psi),\,j_-(w)=j_-(w+\psi),
\]
from \eqref{3.8}, \eqref{3.16} we get
\[
\begin{aligned}
&a_{01}(w,\psi)+a_{11}(w,w,\psi)\geq\langle F_1,\psi\rangle,\\&
a_{01}(w,-\psi)+a_{11}(w,w,-\psi)\geq\langle F_1,-\psi\rangle\quad
\forall \psi\in \mathbf{V}_{8-11}(\Omega),
\end{aligned}
\]
which imply
\begin{equation}\label{3.18}
a_{01}(w,\psi)+a_{11}(w,w,\psi)=\langle F_1,\psi\rangle \quad
\forall \psi\in \mathbf{V}_{8-11}(\Omega).
\end{equation}

When $u\in \mathbf{V}_{10-11}(\Omega)\equiv\{u\in \mathbf{V}(\Omega)
: u|_{\Gamma_{10}\cup\Gamma_{11}}=0\}\left(\subset
K(\Omega)\right)$, the set $\{(u|_{\Gamma_8},u_n|_{\Gamma_9})\}$ is
a subspace of $\textbf{L}^2_\tau(\Gamma_8)\times L^2(\Gamma_9)$,
where  $u_n|_{\Gamma_9}$ is $u|_{\Gamma_9}\cdot n$.

Define a functional $\sigma^*$ on the
set by
\begin{equation}\label{3.19}
\left<\sigma^*,(u|_{\Gamma_8},u_n|_{\Gamma_9})\right>=a_{01}(w,u)+a_{11}(w,w,u)-\langle
F_1,u\rangle\quad \forall u\in \mathbf{V}_{10-11}(\Omega).
\end{equation}
This functional is well defined. Because if $u,u_1\in
\mathbf{V}_{10-11}(\Omega)$ are such that
$(u|_{\Gamma_8},u|_{\Gamma_9})=(u_1|_{\Gamma_8},u_1|_{\Gamma_9})$,
then since $u-u_1\in \mathbf{V}_{8-11}(\Omega)$, by \eqref{3.18}
\[
\begin{aligned}
&a_{01}(w,u-u_1)+a_{11}(w,w,u-u_1)-\langle F_1,u-u_1\rangle=0,
\end{aligned}
\]
that is,
\[
\begin{aligned}
&a_{01}(w,u)+a_{11}(w,w,u)-\langle F_1,u\rangle=a_{01}(w,u_1)+a_{11}(w,w,u_1)-\langle F_1,u_1\rangle,
\end{aligned}
\]
and so  by \eqref{3.19}
\[
\left<\sigma^*,(u|_{\Gamma_8},u_n|_{\Gamma_9})\right>=\left<\sigma^*,(u_1|_{\Gamma_8},u_{1n}|_{\Gamma_9})\right>.
\]
This functional is linear.

Putting $u=w+\psi$, where $\psi\in\mathbf{V}_{10-11}(\Omega)$, and
taking into account
$$j_+(w+\psi)=j_+(w),\quad j_-(w+\psi)=j_-(w),
$$
from \eqref{3.19}, \eqref{3.16} we have
\begin{equation}\label{3.20}
\begin{aligned}
-\left<\sigma^*,(\psi|_{\Gamma_8},\psi_n|_{\Gamma_9})\right>&
=-\left[a_{01}(w,\psi)+a_{11}(w,w,\psi)-\langle
F_1,\psi\rangle\right]
\\&\leq J(w+\psi)-J(w)\\&=
j_\tau(w+\psi)-j_\tau(w)+j_n(w+\psi)-j_n(w)\\
&\leq\int_{\Gamma_8} g_\tau |\psi|_{\Gamma_8}\,dx+\int_{\Gamma_9}
g_n |\psi|_{\Gamma_9}\,dx \quad \forall \psi\in
\mathbf{V}_{10-11}(\Omega).
\end{aligned}
\end{equation}
Putting $u=w-\psi$, in the same way we have
\begin{equation}\label{3.21}
\begin{aligned}
\left<\sigma^*,(\psi|_{\Gamma_8},\psi_n|_{\Gamma_9})\right>&
=\left[a_{01}(w,\psi)+a_{11}(w,w,\psi)-\langle
F_1,\psi\rangle\right]\\&\leq
j_\tau(w-\psi)-j_\tau(w)+j_n(w-\psi)-j_n(w)\\
&\leq\int_{\Gamma_8} g_\tau |\psi|_{\Gamma_8}\,dx+\int_{\Gamma_9}
g_n |\psi|_{\Gamma_9}\,dx \quad \forall \psi\in
\mathbf{V}_{10-11}(\Omega).
\end{aligned}
\end{equation}
By \eqref{3.20}, \eqref{3.21}, we can know that $\sigma^*$ is a
bounded linear functional with a norm not greater than $1$ on a
subspace of $\textbf{L}_{g_\tau}^1(\Gamma_8)\times
L_{g_n}^1(\Gamma_9)$, where $\textbf{L}_{g_\tau}^1(\Gamma_8)$,
$L_{g_n}^1(\Gamma_9)$ are, respectively, the spaces of functions
integrable with weights $g_\tau, g_n$ on $\Gamma_8$ and $\Gamma_9$.
 By the Hahn-Banach theorem the functional is extended as a functional on
$\textbf{L}_{g_\tau}^1(\Gamma_8)\times L_{g_n}^1(\Gamma_9)$ norms of
which is not greater than $1$. Therefore, there exist the elements
$\sigma_\tau\in \textbf{L}_{\f{1}{g_\tau}}^\infty(\Gamma_8)$,
$\|\sigma_\tau\|_{\textbf{L}_{\f{1}{g_\tau}}^\infty(\Gamma_8)}\leq
1$ and $\sigma_n\in L_{\f{1}{g_n}}^\infty(\Gamma_9)$,
$\|\sigma_n\|_{L_{\f{1}{g_n}}^\infty(\Gamma_9)}\leq 1$, which imply
\begin{equation}\label{3.22}
\begin{aligned}
|\sigma_\tau|\leq g_\tau,\quad |\sigma_n|\leq g_n;
\end{aligned}
\end{equation}
and
\begin{equation}\label{3.23}
\left<\sigma^*,(u|_{\Gamma_8},u_n|_{\Gamma_9})\right>=\left(\sigma_\tau,u|_{\Gamma_8}\right)_{\Gamma_8}+\left(\sigma_n,u_n|_{\Gamma_9}\right)_{\Gamma_9}
 \quad \forall u\in \mathbf{V}_{10-11}(\Omega).
\end{equation}

When $u \in \mathbf{V}(\Omega)$, the set
$\{(u_n|_{\Gamma_{10}},u_n|_{\Gamma_{11}})\}$ is a subspace of
$H^{\f{1}{2}}(\Gamma_{10})\times H^{\f{1}{2}}(\Gamma_{11})$.

Define a functional $\sigma^*_1$ on the set $\mathbf{V}(\Omega)$ by
\begin{equation}\label{3.24}
\begin{aligned}
&\left<\sigma^*_1,(u_n|_{\Gamma_{10}},u_n|_{\Gamma_{11}})\right>=\\&\hspace{1cm}a_{01}(w,u)+a_{11}(w,w,u)-\left(\sigma_\tau,u|_{\Gamma_8}\right)_{\Gamma_8}-\left(\sigma_n,u|_{\Gamma_9}\right)_{\Gamma_9}-\langle
F_1,u\rangle\quad \forall u\in \mathbf{V}(\Omega).
\end{aligned}
\end{equation}
 This functional is also well
defined. Because if $u,u^1\in\mathbf{V}(\Omega)$ are such that
$(u_{\Gamma_{10}},u|_{\Gamma_{11}})=(u^1|_{\Gamma_{10}},u^1|_{\Gamma_{11}})$, then since $u-u^1\in
\mathbf{V}_{10-11}(\Omega),$  by \eqref{3.19}, \eqref{3.23}
\[
\begin{aligned}
&a_{01}(w,u-u^1)+a_{11}(w,w,u-u^1)-\left(\sigma_\tau,(u-u^1)|_{\Gamma_8}\right)_{\Gamma_8}-\left(\sigma_n,(u-u^1)|_{\Gamma_9}\right)_{\Gamma_9}-\langle F_1,u-u^1\rangle\\&
=\left<\sigma^*,((u-u^1)|_{\Gamma_{8}},(u-u^1)|_{\Gamma_{9}})\right>-\left(\sigma_\tau,(u-u^1)|_{\Gamma_8}\right)_{\Gamma_8}-\left(\sigma_n,(u-u^1)|_{\Gamma_9}\right)_{\Gamma_9}=0, \end{aligned}
\]
and so  by \eqref{3.24}
\[
\left<\sigma^*_1,(u_n|_{\Gamma_{10}},u_n|_{\Gamma_{11}})\right>=\left<\sigma^*_1,(u^1_n|_{\Gamma_{10}},u^1_n|_{\Gamma_{11}})\right>.
\]
The functional $\sigma^*_1$ is linear. Let us prove continuity of this functional.

 Let $\overline{u}$ is the function corresponding to $\psi\in C_0^\infty(\Gamma_{10})$ by Lemma \ref{l3.1}.
 Then, by Lemma \ref{l3.1} from \eqref{3.24} we have
\begin{equation}\label{3.25.0}
\begin{aligned}
|\left<\sigma^*_1,(\psi,0)\right>|&\leq C\left[ \|w\|_{\textbf{V}}\|\overline{u}\|_{\textbf{V}}+\|w\|_{\textbf{V}}^2\|\overline{u}\|_{\textbf{V}}+(\|\sigma_\tau\|_{L_2(\Gamma_8)}+\|\sigma_\tau\|_{L_2(\Gamma_9)})\|\overline{u}\|+\|F_1\|_{\textbf{V}^*}\|\overline{u}\|_{\textbf{V}}\right]\\&
\leq C\left[ \|w\|_{\textbf{V}}\|+\|w\|_{\textbf{V}}^2+(\|\sigma_\tau\|_{L_2(\Gamma_8)}+\|\sigma_\tau\|_{L_2(\Gamma_9)})+\|F_1\|_{\textbf{V}^*}\right]\cdot\|\psi\|_{H^{\f{1}{2}}(\Gamma_{10})}.\\&
\end{aligned}
\end{equation}
Also assuming that  $\overline{u}$ is the function corresponding to
$\psi\in C_0^\infty(\Gamma_{11})$ by Lemma \ref{l3.1}, we have
\begin{equation}\label{3.25.0.1}
\begin{aligned}
|\left<\sigma^*_1,(0,\psi)\right>|&\leq C\left[
\|w\|_{\textbf{V}}\|\overline{u}\|_{\textbf{V}}+\|w\|_{\textbf{V}}^2\|\overline{u}\|_{\textbf{V}}+(\|\sigma_\tau\|_{L_2(\Gamma_8)}+\|\sigma_\tau\|_{L_2(\Gamma_9)})\|\overline{u}\|+\|F_1\|_{\textbf{V}^*}\|\overline{u}\|_{\textbf{V}}\right]\\&
\leq C\left[
\|w\|_{\textbf{V}}+\|w\|_{\textbf{V}}^2+(\|\sigma_\tau\|_{L_2(\Gamma_8)}+\|\sigma_\tau\|_{L_2(\Gamma_9)})+\|F_1\|_{\textbf{V}^*}\right]\cdot\|\psi\|_{H^{\f{1}{2}}(\Gamma_{11})}.\\&
 \end{aligned}
\end{equation}
Since $H_0^{1/2}(\Gamma_{i})=H^{\f{1}{2}}(\Gamma_{i}), i=10,11,$ (cf. Theorem 11.1 in \cite{lm}), \eqref{3.25.0} and  \eqref{3.25.0.1} show that the functional $\sigma^*_1$  is continuous on the subspace of $H^{\f{1}{2}}(\Gamma_{10})\times H^{\f{1}{2}}(\Gamma_{11})$ mentioned above. Thus, by the Hahn-Banach theorem the functional is extended as a functional on $H^{\f{1}{2}}(\Gamma_{10})\times H^{\f{1}{2}}(\Gamma_{11})$.

Therefore, there exists
an element $(\sigma_{+n}, \sigma_{-n})\in H^{-1/2}(\Gamma_{10})\times H^{-1/2}(\Gamma_{11})$ such that
\begin{equation}\label{3.25}
\left<\sigma^*_1,(u|_{\Gamma_{10}},u|_{\Gamma_{11}})\right>=\left<\sigma_{+n},u|_{\Gamma_{10}}\right>_{\Gamma_{10}}+\left<\sigma_{-n},u|_{\Gamma_{11}}\right>_{\Gamma_{11}}
 \quad \forall u\in \mathbf{V}(\Omega).
\end{equation}

When $\psi\geq 0$ is such that $\psi \in C^\infty_0(\Gamma_{10})$,
let $\overline{u}\in K(\Omega)$ be the function asserted in Lemma
\ref{l3.1}. Putting $u=w+\overline{u}$, by \eqref{3.16} we have
\begin{equation}\label{3.25.1}
\begin{aligned}
&a_{01}(w,\overline{u})+a_{11}(w,w,\overline{u})+J(w+\overline{u})-J(w)-\langle
F_1,\overline{u}\rangle\geq 0.
\end{aligned}
\end{equation}
On the other hand, by \eqref{3.24}, \eqref{3.25} and property of $\overline{u},$
\[
a_{01}(w,\overline{u})+a_{11}(w,w,\overline{u})-\langle
F_1,\overline{u}\rangle=\left<\sigma_{+n},\psi\right>_{\Gamma_{10}}
\]
and so from  \eqref{3.25.1} we have that
\begin{equation}\label{3.25.2}
\left<\sigma_{+n},\psi\right>_{\Gamma_{10}}+J(w+\overline{u})-J(w)\geq
0.
\end{equation}
By \eqref{3.7},  \eqref{3.8} and property of $\overline{u},$
\[
J(w+\overline{u})-J(w)=\left<g_{+n},\psi\right>_{\Gamma_{10}},
\]
and combining with \eqref{3.25.2} we have
\[
\left<\sigma_{+n},\psi\right>_{\Gamma_{10}}+(g_{+n},\psi)_{\Gamma_{10}}\geq
0,
\]
that is,
\begin{equation}\label{3.26}
\sigma_{+n}+g_{+n}\geq 0.
\end{equation}

When $\psi\leq 0$ is such that $\psi \in C^\infty_0(\Gamma_{11})$,
let $\overline{u}\in K(\Omega)$ be the function asserted in Lemma
\ref{l3.1}. Then, in the same way we have that
\[
\left<\sigma_{-n},-\psi\right>_{\Gamma_{11}}-(g_{-n},-\psi)_{\Gamma_{11}}\geq
0,
\]
that is,
\begin{equation}\label{3.27}
\sigma_{-n}-g_{-n}\leq 0.
\end{equation}

 From  \eqref{3.24}, \eqref{3.25}, we have
\begin{equation}\label{3.28}
\begin{aligned}
a_{01}(w,u)+a_{11}(w,w,u)&-(\sigma_\tau,u_\tau)_{\Gamma_8}-(\sigma_n,u)_{\Gamma_9}-\left<\sigma_{+n},u\right>_{\Gamma_{10}}-\left<\sigma_{-n},u\right>_{\Gamma_{11}}\\&
=\langle F_1,u\rangle\quad  \forall u\in \mathbf{V}(\Omega).
\end{aligned}
\end{equation}

Putting $u=0$ in \eqref{3.16} and taking into account \eqref{3.28}
with $u=w$, we have
\[
\begin{aligned}
(\sigma_\tau,w)_{\Gamma_8}+(\sigma_n,w)_{\Gamma_9}&+\left<\sigma_{+n},w_n\right>_{\Gamma_{10}}+\left<\sigma_{-n},w_n\right>_{\Gamma_{11}}\\&
+j_\tau(w)+j_n(w)+j_+(w)+j_-(w)\leq 0,
\end{aligned}
\]
that is,
\begin{equation}\label{3.29}
\begin{aligned}
\int_{{\Gamma_8}}(\sigma_\tau w_\tau&+g_\tau|w_\tau|)\,
ds+\int_{{\Gamma_9}}(\sigma_n w_n+g_n|w_n|)\, ds\\&
+\left<\sigma_{+n}+g_{+n},w_n\right>_{\Gamma_{10}}+\left<\sigma_{-n}-g_{-n},w_n\right>_{\Gamma_{11}} \leq 0.
\end{aligned}
\end{equation}
Since on $\Gamma_8, \Gamma_9$, $\Gamma_{10}$ and $\Gamma_{11}$,
respectively, $w_\tau=v_\tau$, $w_n=v_n$, $w_n=v_n\geq 0$ and
$w_n=v_n\leq 0$, taking into account \eqref{3.22},
\eqref{3.26},\eqref{3.27}, by  \eqref{3.29} we have
\begin{equation}\label{3.30}
\begin{aligned}
&\sigma_\tau v_\tau+g_\tau|v_\tau|=0,\quad \sigma_n
v_n+g_n|v_n|=0,\\&\left<\sigma_{+n}+g_{+n},v_n\right>=0,\quad\left<\sigma_{-n}-g_{-n},v_n\right>=0.
\end{aligned}
\end{equation}
 Therefore, by virtue of \eqref{3.17}, \eqref{3.22},
\eqref{3.26}-\eqref{3.28}, \eqref{3.30},
 we come to the conclusion. $\square$\vspace{.4cm}

 Taking $(v\cdot\nabla)v=\mbox{rot}\,v\times
v+\frac{1}{2}\mbox{grad}|v|^2$ into account and putting $v=w+U$, by
\eqref{3.1}, \eqref{3.2} and Assumption \ref{a3.1} we can see that
smooth solutions $v$ of problem \eqref{2.7}, \eqref{2.9} satisfy the
following.
\begin{equation}\label{3.31}
\left\{\begin{aligned} &v-U=w\in K(\Omega),\\&
 2\nu(\varepsilon(w),\varepsilon(u))+\langle\mbox{rot}\,w\times
w,u\rangle+\langle\mbox{rot}\,U\times
w,u\rangle+\langle\mbox{rot}\,w\times
U,u\rangle\\&\hspace*{1cm}+2\nu(k(x)w,u)_{\Gamma_2}+2\nu(S\tilde{w},\tilde{u})_{\Gamma_3}
+2(\alpha(x)w,u)_{\Gamma_5}+\nu(k(x)w,u)_{\Gamma_7}\\&
\hspace*{1cm}-2(\varepsilon_{n\tau}(w+U),u)_{\Gamma_8}
+\big(p+\f{1}{2}|v|^2-2\varepsilon_{nn}(w+U),u_n\big)_{\Gamma_9\cup\Gamma_{10}\cup\Gamma_{11}}\\&\hspace*{.5cm}=-2\nu(\varepsilon(U),\varepsilon(u))-\langle\mbox{rot}\,U\times
U,u\rangle-2\nu(k(x)U,u)_{\Gamma_2}-2\nu(S\tilde{U},\tilde{u})_{\Gamma_3}\\
&\hspace*{1cm}-2(\alpha(x)U,u)_{\Gamma_5}-\nu(k(x)U,u)_{\Gamma_7}+\langle
f,u\rangle+\sum_{i=2,4,7}\langle\phi_i,u_n\rangle_{\Gamma_i}+\sum_{i=3,5,6}\langle\phi_i,u\rangle_{\Gamma_i}\\&\hspace{5cm}
 \forall u\in\textbf{ V}(\Omega),
\\&|\sigma_\tau^t(v)|\leq g_\tau,\,
\sigma_\tau^t\cdot v_\tau+g_\tau |v_\tau|=0\quad \mbox{on}\,\,
\Gamma_8,\\
&|\sigma_n^t(v)|\leq g_n,\, \sigma_n^t(v) v_n+g_n |v_n|=0\quad
\mbox{on}\,\, \Gamma_9
,\\
 &
 \sigma^t_n(v)+g_{+n}\geq 0,\quad
(\sigma^t_n(v)+g_{+n}) v_n=0\quad \mbox{on}\,\, \Gamma_{10},
\\
&\sigma^t_n(v)-g_{-n}\leq 0,\quad (\sigma^t_n(v)-g_{-n}) v_n=0\quad
\mbox{on}\,\, \Gamma_{11}.
\end{aligned}
\right.
\end{equation}
Define $a_{02}(\cdot, \cdot), a_{12}(\cdot,\cdot,\cdot)$ and $F_2\in
V^*$ by
\begin{equation}\label{3.32}
\begin{aligned}
&
a_{02}(w,u)=2\nu(\varepsilon(w),\varepsilon(u))+\langle\mbox{rot}\,U\times
w,u\rangle+\langle\mbox{rot}\,w\times
U,u\rangle+2\nu(k(x)w,u)_{\Gamma_2}\\&
\hspace{2.1cm}+2\nu(S\tilde{w},\tilde{u})_{\Gamma_3}+2(\alpha(x)w,u)_{\Gamma_5}+\nu(k(x)w,u)_{\Gamma_7}\quad
 \forall w,u\in\textbf{ V}(\Omega),
\\&a_{12}(w,u,v)=\langle\mbox{rot}\,w\times
u,v\rangle\quad \forall w,u,v\in \mathbf{V}(\Omega),\\& \langle
F_2,u\rangle=-2\nu(\varepsilon(U),\varepsilon(u))-\langle\mbox{rot}\,U\times
U,u\rangle-2\nu(k(x)U,u)_{\Gamma_2}-2\nu(S\tilde{U},\tilde{u})_{\Gamma_3}-2(\alpha(x)U,u)_{\Gamma_5}\\
&\hspace*{1cm}-\nu(k(x)U,u)_{\Gamma_7}+\langle
f,u\rangle+\sum_{i=2,4,7}\langle\phi_i,u_n\rangle_{\Gamma_i}+\sum_{i=3,5,6}\langle\phi_i,u\rangle_{\Gamma_i}\quad
\forall u\in \mathbf{V}(\Omega).
\end{aligned}
\end{equation}

Then, taking into account
\[
\begin{aligned}
&\sigma^t_\tau(v)=2\nu\varepsilon_{n\tau}(v),\quad
\sigma^t_n(v)=-(p+\f{1}{2}|v|^2)+2\nu\varepsilon_{nn}(v)
\end{aligned}
\]
and \eqref{3.31}, we introduce the following variational formulation
for problem \eqref{2.7}, \eqref{2.9}.\vspace{.2cm}

\textbf{Problem
II-VE.} Find $(v, \sigma^t_\tau, \sigma^t_n, \sigma^t_{+n},
\sigma^t_{-n})\in \big(U+K(\Omega)\big)\times
\textbf{L}^2_\tau(\Gamma_8)\times L^2(\Gamma_9)\times
H^{-\f{1}{2}}(\Gamma_{10})\times H^{-\f{1}{2}}(\Gamma_{11})$ such that
\begin{equation}\label{3.33}
\left\{\begin{aligned} &v-U=w\in K(\Omega),\\&
a_{02}(w,u)+a_{12}(w,w,u)-(\sigma^t_\tau,u_\tau)_{\Gamma_8}-(\sigma^t_n,u_n)_{\Gamma_9}\\&
\hspace{1.4cm}-\left<\sigma^t_{+n},u_n\right>_{\Gamma_{10}}-\left<\sigma^t_{-n},u_n\right>_{\Gamma_{11}}=\langle
F_2,u\rangle\quad  \forall u\in \textbf{V}(\Omega),
\\
&|\sigma^t_\tau|\leq g_\tau,\, \sigma^t_\tau\cdot v_\tau+g_\tau
|v_\tau|=0\quad \mbox{on}\,\, \Gamma_8,\\
&|\sigma^t_n|\leq g_n,\, \sigma^t_n v_n+g_n |v_n|=0\quad
\mbox{on}\,\, \Gamma_9,
\\&\sigma^t_{+n}+g_{+n}\geq 0,\quad
\left<\sigma^t_{+n}+g_{+n}, v_n\right>_{\Gamma_{10}}=0\quad \mbox{on}\,\, \Gamma_{10},
\\
&\sigma^t_{-n}-g_{-n}\leq 0,\quad \left<\sigma^t_{-n}-g_{-n}, v_n\right>_{\Gamma_{11}} =0\quad
\mbox{on}\,\, \Gamma_{11}.
\end{aligned}
\right.
\end{equation}

In the same way as Theorem \ref{t3.0} we have
\begin{theorem}\label{t3.2.0} Assume 1), 2) of Assumption \ref{a3.1}. If a
solution smooth enough ($v\in\textbf{H}^2(\Omega), f\in
\textbf{L}^2(\Omega)$), then Problem II-VE is equivalent to problem
\eqref{2.7}, \eqref{2.9}. In addition, if among
$\Gamma_i,i=2,4,6,7,9-11,$ at least one is nonempty, then $p$ of
problem \eqref{2.7}, \eqref{2.8} is unique.
\end{theorem}

Then, in the same way as Problem I we get Problem II-VI formulated
by a variational inequality and can prove that the problem is
equivalent to Problem II-VE.\vspace{.1cm}

\textbf{Problem II-VI.} Find $v=w+U$ such that
\begin{equation}\label{3.34}
\begin{aligned}
&a_{02}(w,u-w)+a_{12}(w,w,u-w)+J(u)-J(w)\geq\langle
F_2,u-w\rangle\quad  \forall u\in \textbf{v}(\Omega),
\end{aligned}
\end{equation}
where $a_{02}, a_{12}, F_2$ are in \eqref{3.32} and $J$ is defined
by \eqref{3.7}, \eqref{3.8}. \vspace*{.3cm}

\begin{theorem}\label{t3.2} If $(v, \sigma^t_\tau, \sigma^t_n, \sigma^t_{+n},
\sigma^t_{-n})$ is a solution to Problem II-VE, then $v$ is a
solution to Problem II-VI. Inversely, if $v$ is a solution to
Problem II-VI, then there exist $\sigma^t_\tau, \sigma^t_n,
\sigma^t_{+n}, \sigma^t_{-n}$ such that $(v, \sigma^t_\tau,
\sigma^t_n, \sigma^t_{+n}, \sigma^t_{-n})$  is a solution to Problem
I-VE.
\end{theorem}

\begin{remark}\label{r3.2}
 Boundary
condition $\nu\f{\p v}{\p n}-pn=0$ often called ``do nothing" or
``free outflow" boundary condition, results from variational
principle and does not have a real physical meaning but is rather
used in truncating large physical domains to smaller computational
domains by assuming parallel flow (cf. \cite{b}). The condition (7)
in \eqref{2.8}(corresponding (7) in \eqref{2.9}) is rather different
from ``do nothing" condition. Assuming that the
flow is orthogonal on $\Gamma_7$ and applying Theorem \ref{t2.2}, we get a variational formulation,
 and so to convert from the variational formulation to the original problem we use such a
 condition. (For more
detail refer to Remark 2.1 in \cite{kc3}.) If the flow, in addition,
is parallel in a near the boundary, then condition (7) in
\eqref{2.8} is same with ``do nothing" condition. In point of view
of pure mathematics, to reflect correctly ``do nothing" condition in
variational formulation we can use other variational formulation
assuming $\Gamma_6=\emptyset$. Bellow we show that.
\end{remark}\vspace*{.2cm}

 Now, we consider the cases that $\Gamma_6=\varnothing$ and for
convenience $h_i=0,\,\,i=4, 5, 8, 9,$ in \eqref{2.8}. Let
$\mathbf{V}_{\Gamma7}(\Omega)=\{u\in
\mathbf{H}^1(\Omega):\mbox{div}\,u=0,\, u|_{\Gamma_1}=0,\,
u_\tau|_{(\Gamma_2\cup\Gamma_4\cup\Gamma_9\cup\Gamma_{10}\cup\Gamma_{11})}=0,\,
u\cdot n|_{(\Gamma_3\cup\Gamma_5\cup\Gamma_8)}=0\}$ and
$\mathbf{V}_{\Gamma17}(\Omega)=\{u\in
\mathbf{H}^1(\Omega):\mbox{div}\,u=0,\,
u_\tau|_{(\Gamma_2\cup\Gamma_4\cup\Gamma_9\cup\Gamma_{10}\cup\Gamma_{11})}=0,\,
u\cdot n|_{(\Gamma_3\cup\Gamma_5\cup\Gamma_8)}=0\}$. By Theorem
\ref{2.1} and \ref{2.2} for $v\in
\mathbf{H}^2(\Omega)\cap\mathbf{V}_{\Gamma17}(\Omega)$ and $u\in
\mathbf{V}_{\Gamma7}(\Omega)$
\begin{equation}\label{3.35}
\begin{aligned}
-(\Delta& v, u)=(\nabla v, \nabla u)-\left(\f{\p v}{\p
n},u\right)_{\cup_{i=2}^{11}\Gamma_i}\\& =(\nabla v, \nabla
u)+(k(x)v,u)_{\Gamma_2}-(\mbox{rot}\,v\times n,
u)_{\Gamma_3}+(S\tilde{v},\tilde{u})_{\Gamma_3}-(\varepsilon_{nn}(v),u\cdot
n)_{\Gamma_4} -2(\varepsilon_{n\tau}(v),u)_{\Gamma_5}\\&
\quad-(S\tilde{v},\tilde{u})|_{\Gamma_5}-\left(\f{\p v}{\p
n},u\right)_{\Gamma_7}-2(\varepsilon_{n\tau}(w),u)_{\Gamma_8}-(S\tilde{v},\tilde{u})_{\Gamma_8}-(\varepsilon_{nn}(v),u_n)_{\Gamma_9\cup\Gamma_{10}\cup\Gamma_{11}}.
\end{aligned}
\end{equation}

 Using \eqref{3.35}, \eqref{3.2} we get a variational
formulation for problem \eqref{2.7}, \eqref{2.8} with $(-p\cdot
n+\nu\f{\p v}{\p n})|_{\Gamma_7}=\phi_7\in
\mathbf{H}^{-\f{1}{2}}(\Gamma_7)$ instead of the condition
 (7) of \eqref{2.8}:\vspace*{.3cm}

\textbf{Problem I'-VE.} Find $(v, \sigma_\tau, \sigma_n,
\sigma_{+n}, \sigma_{-n})\in \big(U+K(\Omega)\big)\times
\textbf{L}^2_\tau(\Gamma_8)\times L^2(\Gamma_9)\times
H^{-\f{1}{2}}(\Gamma_{10})\times H^{-\f{1}{2}}(\Gamma_{11})$ such that
\begin{equation}\label{3.36}
\left\{\begin{aligned} & v|_{\Gamma_1}=h_1,\\&
 \nu(\nabla v, \nabla u)+(v\cdot \nabla)v+\nu(k(x)v,u)_{\Gamma_2}+\nu(S\tilde{v},\tilde{u})_{\Gamma_3}+2(\alpha(x)v,u)_{\Gamma_5}
 -\nu(S\tilde{v},\tilde{u})_{\Gamma_5}\\&\hspace*{1.5cm}-\nu(S\tilde{v},\tilde{u})_{\Gamma_8}-(\sigma_\tau,u_\tau)_{\Gamma_8}-(\sigma_n,u_n)_{\Gamma_9}
 -\left<\sigma_{+n},u_n\right>_{\Gamma_{10}}-\left<\sigma_{-n},u_n\right>_{\Gamma_{11}}\\&
 \hspace*{1cm}=\langle f,u\rangle+\sum_{i=2,4}\langle \phi_i,u_n\rangle_{\Gamma_i}
 +\sum_{i=3,5,7}\langle \phi_i,u\rangle_{\Gamma_i}\quad \forall u\in
\textbf{V}(\Omega),\\& |\sigma_\tau|\leq g_\tau,\, \sigma_\tau\cdot
v_\tau+g_\tau
|v_\tau|=0\quad \mbox{on}\,\, \Gamma_8,\\
&|\sigma_n|\leq g_n,\, \sigma_n v_n+g_n |v_n|=0\quad \mbox{on}\,\,
\Gamma_9,
\\&\sigma_{+n}+g_{+n}\geq 0,\quad
\left<\sigma_{+n}+g_{+n}, v_n\right>_{\Gamma_{10}}=0\quad \mbox{on}\,\, \Gamma_{10},
\\
&\sigma_{-n}-g_{-n}\leq 0,\quad \left<\sigma_{-n}-g_{-n}, v_n\right>_{\Gamma_{11}}=0\quad
\mbox{on}\,\, \Gamma_{11}.
\end{aligned}
\right.
\end{equation}
\vspace*{.3cm}

In the same way as Problem I we get the below equivalent
formulations of Problem III for the Stokes equation with boundary
condition \eqref{2.8}.

\textbf{Problem III-VE.} Find $(v, \sigma_\tau, \sigma_n,
\sigma_{+n}, \sigma_{-n})\in \big(U+K(\Omega)\big)\times
\textbf{L}^2_\tau(\Gamma_8)\times L^2(\Gamma_9)\times
H^{-\f{1}{2}}(\Gamma_{10})\times H^{-\f{1}{2}}(\Gamma_{11})$ such that
\begin{equation}\label{3.37}
\left\{\begin{aligned} &v-U=w\in K(\Omega),\\&
a_{03}(w,u)-(\sigma_\tau,u_\tau)_{\Gamma_8}-(\sigma_n,u_n)_{\Gamma_9}-\left<\sigma_{+n},u_n\right>_{\Gamma_{10}}-\left<\sigma_{-n},u_n\right>_{\Gamma_{11}}=\langle
F_3,u\rangle\quad  \forall u\in \textbf{V}(\Omega),
\\&
|\sigma_\tau|\leq g_\tau,\, \sigma_\tau\cdot v_\tau+g_\tau
|v_\tau|=0\quad \mbox{on}\,\, \Gamma_8,\\
&|\sigma_n|\leq g_n,\, \sigma_n v_n+g_n |v_n|=0\quad \mbox{on}\,\,
\Gamma_9,
\\&\sigma_{+n}+g_{+n}\geq 0,\quad
\left<\sigma_{+n}+g_{+n}, v_n\right>=0\quad \mbox{on}\,\, \Gamma_{10},
\\
&\sigma_{-n}-g_{-n}\leq 0,\quad \left<\sigma_{-n}-g_{-n}, v_n\right>=0\quad
\mbox{on}\,\, \Gamma_{11},
\end{aligned}
\right.
\end{equation}
where
\begin{equation}\label{3.38}
\begin{aligned}
&
a_{03}(w,u)=2\nu(\varepsilon(w),\varepsilon(u))+2\nu(k(x)w,u)_{\Gamma_2}\\&
\hspace{1.6cm}+2\nu(S\tilde{w},\tilde{u})_{\Gamma_3}+2(\alpha(x)w,u)_{\Gamma_5}+\nu(k(x)w,u)_{\Gamma_7}\quad
\forall w,u\in \mathbf{V}(\Omega),\\
& \langle
F_3,u\rangle=-2\nu(\varepsilon(U),\varepsilon(u))-2\nu(k(x)U,u)_{\Gamma_2}-2\nu(S\tilde{U},\tilde{u})_{\Gamma_3}
-2(\alpha(x)U,u)_{\Gamma_5}\\& \hspace*{1.3cm}
-\nu(k(x)U,u)_{\Gamma_7}+\langle
f,u\rangle+\sum_{i=2,4,7}\langle\phi_i,u_n\rangle_{\Gamma_i}+\sum_{i=3,5,6}\langle\phi_i,u\rangle_{\Gamma_i}\quad
\forall u\in \mathbf{V}(\Omega).
\end{aligned}
\end{equation}

 \textbf{Problem III-VI} Find $v$ such that
\begin{equation}\label{3.39}
\begin{aligned}
& v-U=w\in K(\Omega),\\& a_{03}(w,u-w)+J(u)-J(w)\geq\langle
F_3,u-w\rangle\quad \forall u\in \mathbf{V}(\Omega),
\end{aligned}
\end{equation}
where the functionals $J$ is defined by \eqref{3.7},
\eqref{3.8}.\vspace*{.3cm}

\section{Existence, uniqueness and estimates of solutions to variational inequalities}
\setcounter{equation}{0}

In this section we study some variational inequalities for the
problems in Section 3.

\begin{theorem}\label{t4.1} Let $X, X_1$ be real separable Hilbert spaces such that $X\hookrightarrow\hookrightarrow X_1$, and $X^*$ be dual space of $X$. Assume the
followings.

1) $J\in(X\rightarrow [0, +\infty])$ is a proper lower
semi-continuous convex functional such that $J(0_X)=0$.

2)  $a_0(\cdot,\cdot)\in (X \times X\rightarrow R)$ is a bilinear
form such that
\[
\begin{aligned}
 & |a_0(u,v)|\leq K \|u\|_X\|v\|_X\quad \forall
 u,v\in X,\\&
 |a_0(u,u)|\geq \alpha\|u\|_X^2\quad  \exists\alpha>0, \forall
 u\in X.\\
\end{aligned}
\]

3)  $a_1(\cdot,\cdot,\cdot)\in (X_1\times X\times X\rightarrow R)$
is a triple linear functional such that
\[
\begin{aligned}
 & a_1(w,u,u)= 0\quad \forall
 w\in X_1, \forall u\in X,\\&
 |a_1(w,u,v)|\leq K\|w\|_{X_1}\|u\|_X\|v\|_X,\quad \forall
 w\in X_1, \forall u,v\in X.
 \end{aligned}
 \]

 Then for $f\in X^*$ there exists a solution to the variational inequality
 \begin{equation}\label{4.1}
a_0(v,u-v)+a_1(v,v,u-v)+J(u)-J(v)\geq \langle f, u-v\rangle\quad
\forall u\in X
\end{equation}
and all solutions $v$ satisfy the estimate
\begin{equation}\label{4.2}
\|v\|_X\leq \f{1}{\alpha}\|f\|_{X^*}.
\end{equation}
In addition to, if
\begin{equation}\label{4.3}
\f{Kc}{\alpha^2}\|f\|_{X^*}<1,
\end{equation}
then solution is unique, where $c$ is a constant in
$\|\cdot\|_{X_1}\leq c\|\cdot\|_X$.
\end{theorem}
$\mathbf{Proof.}$\, Fixing $w\in {X_1}$, let us consider a
variational inequality
\begin{equation}\label{4.4}
a_0(v,u-v)+a_1(w,v,u-v)+J(u)-J(v)\geq \langle f, u-v\rangle\quad
\forall u\in X,
\end{equation}
where $f\in X^*$. There exists a unique solution to \eqref{4.4} (cf.
Theorem 10.5 in \cite{bc}). Let $v_1, v_2$ be the solutions
corresponding to $f_1, f_2$ instead of $f$. Then, under
consideration of condition 2) it is easy to verify that
\begin{equation}\label{4.5}
\|v_1-v_2\|_X\leq \f{1}{\alpha}\|f_1-f_2\|_{X^*}.
\end{equation}

Now, let us consider the operator which maps $w$ to the solution $v$
of \eqref{4.4}
\[
 T\in(X_1\rightarrow X):
w\rightarrow T(w)=v.
\]
Taking into account condition 1), we can easily verify that the
solution corresponding to $f=0_{X^*}$ is $0_X$. Thus, from
\eqref{4.5} we have
\begin{equation}\label{4.6}
\|v\|_X\leq \f{1}{\alpha}\|f\|_{X^*}\quad \forall w\in {X_1}.
\end{equation}
Note that this estimate is independent from $w$.

 Denote by $v_1$ and $v_2$,
respectively, the solutions to \eqref{4.4} corresponding to $w_1$
and $w_2$. Then
\begin{equation}\label{4.7}
\begin{aligned}
&a_0(v_1,u-v_1)+a_1(w_1,v_1,u-v_1)+J(u)-J(v_1)\geq \langle f,
u-v_1\rangle\quad
\forall u\in X,\\
&a_0(v_2,u-v_2)+a_1(w_2,v_2,u-v_2)+J(u)-J(v_2)\geq \langle f,
u-v_2\rangle\quad \forall u\in X.
\end{aligned}
\end{equation}
Putting $u=v_2$ and $u=v_1$, respectively, in the first formula and
the second one of \eqref{4.7}, and adding two formulae, we get
\begin{equation}\label{4.8}
a_0(v_1-v_2,v_2-v_1)+a_1(w_1,v_1,v_2-v_1)+a_1(w_2,v_2,v_1-v_2)\geq
0.
\end{equation}
From \eqref{4.8}, the conditions 2), 3) of Theorem and \eqref{4.6},
we get
\[
\begin{aligned}
\|v_2-v_1\|_X^2&\leq\f{1}{\alpha}|a_1(w_1,v_1,v_2-v_1)-a_1(w_2,v_1,v_2-v_1)\\&\hspace{2cm}+a_1(w_2,v_1,v_2-v_1)-a_1(w_2,v_2,v_2-v_1)|\\
&\leq\f{1}{\alpha}|a_1(w_1-w_2,v_1,v_2-v_1)|+\f{1}{\alpha}|a_1(w_2,v_2-v_1,v_2-v_1)|\\&
\leq
\f{K}{\alpha}\|w_1-w_2\|_{X_1}\|v_1\|_X\|v_2-v_1\|_X\\
&\leq \f{K\|f\|_{X^*}}{\alpha^2}\|w_1-w_2\|_{X_1}\|v_2-v_1\|_X\quad
\forall
 w_1,w_2\in X_1,
\end{aligned}
\]
which implies
\begin{equation}\label{4.9}
\|v_2-v_1\|_X\leq \f{K\|f\|_{X^*}}{\alpha^2}\|w_1-w_2\|_{X_1}\quad
\forall
 w_1,w_2\in X_1.
\end{equation}
By \eqref{4.6}, \eqref{4.9} and Schauder fixed-point theorem(cf.
Theorem 2.A in \cite{z}) there exists a solution to \eqref{4.1}. And
any solution is a fixed point of operator $T,$ and by \eqref{4.6}
all solutions satisfy the estimate \eqref{4.2}.

If \eqref{4.3} holds, then the operator $T:w\in X \rightarrow v\in
X$ is contract, and so we come to the last conclusion.
$\square$\vspace*{.3cm}

Let us study variational inequalities when the condition 3) of the
above theorem is weakened.

\begin{theorem}\label{t4.2} Let $X$ be a real separable Hilbert space. Assume the
followings.

1) Condition 1) of Theorem \ref{t4.1} holds.

2) Condition 2) of Theorem \ref{t4.1} holds.

3)  $a_1(\cdot,\cdot,\cdot)\in (X\times X\times X\rightarrow R)$ is
a triple linear functional such that
\[
  |a_1(w,u,v)|\leq K\|w\|_X\|u\|_X\|v\|_X,\quad \forall
 w,u,v\in X.
 \]

 If $f$ is small enough, then in $\mathscr{O}_M(0_X)$, where $M$ is determined in \eqref{4.18}, there exists a unique solution to the variational inequality
\begin{equation}\label{4.10}
a_0(v,u-v)+a_1(v,v,u-v)+J(u)-J(v)\geq \langle f, u-v\rangle\quad
\forall u\in X.
\end{equation}
\end{theorem}

$\mathbf{Proof.}$\, Fixing $w\in X$, let us consider a variational
inequality
\begin{equation}\label{4.11}
a_0(v,u-v)+a_1(w,w,u-v)+J(u)-J(v)\geq \langle f, u-v\rangle\quad
\forall u\in X,
\end{equation}
where $f\in X^*$. Defining an element $a_1(w)\in X^*$ by
\[
\langle a_1(w),u\rangle=a_1(w,w,u)\quad \forall u\in X,
\]
by condition 3) we have
\begin{equation}\label{4.12}
\|a_1(w)\|_{X^*}\leq K\|w\|_X^2\quad \forall w\in X,
\end{equation}
 Then, \eqref{4.11} is rewritten as follows.
\begin{equation}\label{4.13}
a_0(v,u-v)+J(u)-J(v)\geq \langle f-a_1(w), u-v\rangle\quad \forall
u\in X.
\end{equation}
By the same argument as Theorem \ref{t4.1}, there exists a unique
solution $v_w$ to \eqref{4.13} and
\begin{equation}\label{4.14}
\|v_w\|\leq \f{1}{\alpha}(\|f\|_{X^*}+\|a_1(w)\|_{X^*})\leq
\f{1}{\alpha}(\|f\|_{X^*}+K\|w\|_X^2),
\end{equation}
where \eqref{4.12} was used.

Now, let us consider the operator which maps $w$ to the solution of
\eqref{4.13}
\[
 T\in(X\rightarrow X):
w\rightarrow T(w)=v
\]
Denote by $v_1$ and $v_2$, respectively, the solutions to
\eqref{4.11} corresponding to $w_1,w_2\in \mathscr{O}_M(0_X)$, where
$M$ is determined bellow. Then
\begin{equation}\label{4.15}
\|v_1-v_2\|_X\leq \f{1}{\alpha}\|a_1(w_1)-a_1(w_2)\|_{X^*}.
\end{equation}
By condition 3)
\begin{equation}\label{4.16}
\begin{aligned}
 \|a_1(w_1)-a_1(w_2)\|_{X^*}&\leq
 K\left(\|w_2-w_1\|_X\|w_2\|_X+\|w_1\|_X\|w_1-w_2\|_X\right).
\end{aligned}
\end{equation}
Thus, by \eqref{4.15}, \eqref{4.16}
\begin{equation}\label{4.17}
\begin{aligned}
 \|v_1-v_2\|_X&\leq
\f{K}{\alpha}\left(\|w_2-w_1\|_X\|w_2\|_X+\|w_1\|_X\|w_1-w_2\|_X\right)\\
&\leq\f{2KM}{\alpha}\|w_2-w_1\|_X\quad \forall w_1,w_2\in
\mathscr{O}_M(0_X).
\end{aligned}
\end{equation}
Therefore, if $M$ is taken satisfied (If $\alpha$ is large and
$\|f\|_{X^*}$ is small enough, then such choosing is possible.)
\begin{equation}\label{4.18}
\left\{\begin{aligned}
&M=\f{1}{\alpha}(\|f\|_{X^*}+KM^2),\\
&\f{2KM}{\alpha}<1,
\end{aligned}\right.
\end{equation}
 then by \eqref{4.14}, \eqref{4.17} the operator $T$ on
$\mathscr{O}_M(0_X)$ is contract, and so there exists a unique
solution to \eqref{4.10}. $\square$\vspace*{.3cm}

\begin{theorem}\label{t4.3} Let $X$ be a real separable Hilbert space and $X^*$ be its dual space. Assume that

 1) $J\in(X\rightarrow R)$ is a finite weak continuous convex functional, $J_\varepsilon\in (X\rightarrow R)$
  is convex such that
\[
\begin{aligned}
 & J_\varepsilon(v)\rightarrow  J(v)\quad \mbox{uniformly on}\,\, X\,\, \mbox{as}\,\, \varepsilon\rightarrow
 0,\\
 & \mbox{Gateaux derivative}\,\, DJ_\varepsilon\equiv A_\varepsilon\in(X\rightarrow X^*) \,\, \mbox{is weak
 continuous and }\,\,A_\varepsilon(0_X)=0_{X^*};
\end{aligned}
\]

2)  $a(\cdot,\cdot, \cdot)\in (X\times X \times X\rightarrow R)$ is
a form such that
\[
\begin{aligned}
 & \mbox{when}\,\, w\in X, (u,v)\rightarrow a(w;u,v)\,\, \mbox{is bilinear on}\,\, X\times
 X,\\
&a(v,v,v)\geq \alpha\|v\|^2_X\quad \exists \alpha>0, \forall
 v\in X\,\,\text{and}\\&
 \mbox{when}\,\,v_m \rightharpoonup v\,\, \mbox{weakly in}\,\, X,\,\,a(v_m,v_m,u)\rightarrow a(v,v,u) \,\,\forall u\in
 X\,\,\mbox{and}\\&\hspace{4.2cm}   \liminf_{m\rightarrow \infty} a(v_m,v_m,v_m)\geq a(v,v,v).
 \end{aligned}
 \]
 Then for $f\in X^*$ there exists a solution to a variational inequality
 \begin{equation}\label{4.19}
a(v,v,u-v)+J(u)-J(v)\geq \langle f, u-v\rangle\quad \forall u\in X
\end{equation}
satisfying an estimate
\begin{equation}\label{4.20}
\|v\|_X\leq \f{1}{\alpha}\|f\|_{X^*}.
\end{equation}
\end{theorem}
$\mathbf{Proof.}$\, First let us prove existence of a solution to a
variational equation
\begin{equation}\label{4.21}
a(v,v,u)+\langle A_\varepsilon(v),u\rangle=\langle f, u\rangle\quad
\forall u\in X.
\end{equation}

We will do it as Theorem 1.2 in ch. 4 of \cite{gr}. Let $\{w_n\}$ be
a base of $X$ and denote by $X_m$ the subspace of $X$ spanned by
$w_1, \cdots, w_m$.

We find $v_m=\sum_{i=1}^m\nu_iw_i\in X_m$ satisfying
\begin{equation}\label{4.22}
a(v_m,v_m,u)+\langle A_\varepsilon(v_m),u\rangle=\langle f,
u\rangle\quad \forall u\in X_m.
\end{equation}
Define $\Phi_m\in(X_m\rightarrow X_m)$ by
\begin{equation}\label{4.23}
(\Phi_m(v),w_i)=a(v,v,w_i)+\langle
A_\varepsilon(v),w_i\rangle-\langle f, w_i\rangle,\quad 1\leq i\leq
m.
\end{equation}
Since Gateaux derivative of convex functional is monotone (cf. Lemma
4.10, ch. 3 in \cite{ggz}) and $A_\varepsilon(0_X)=0_{X^*}$,
\[
\langle A_\varepsilon(u)-A_\varepsilon(0_X),u-0_X\rangle=\langle
A_\varepsilon(u),u\rangle\geq 0\quad \forall u\in X.
\]
 Thus,
\begin{equation}\label{4.24}
a(u,u,u)+\langle A_\varepsilon(u),u\rangle\geq \alpha\|u\|^2_X\quad
\forall u\in X.
\end{equation}
From \eqref{4.23}, \eqref{4.24} we get
\begin{equation}\label{4.25}
(\Phi_m(v),v)\geq (\alpha\|v\|_X-\|f\|_{X^*})\|v\|_X\quad \forall
v\in X_m.
\end{equation}
Therefore,
\[
(\Phi_m(v),v)\geq 0\quad \forall v\in X \,\,\mbox{with}\,\,
\|v\|_X=\f{\|f\|_{X^*}}{\alpha}.
\]
And $\Phi_m$ is continuous in $X_m$ by virtue of the assumption 2).
Thus, there exists a solution $v_{\varepsilon m}$ to problem
\eqref{4.22}. By \eqref{4.25} for all solution $v_{\varepsilon m}$
to \eqref{4.22}
\[
0=(\Phi_m(v_{\varepsilon m}),v_{\varepsilon m})\geq
(\alpha\|v_{\varepsilon m}\|_X-\|f\|_{X^*})\|v_{\varepsilon m}\|_X,
\]
which implies
\begin{equation}\label{4.26}
\|v_{\varepsilon m}\|_X\leq \f{1}{\alpha}\|f\|_{X^*}.
\end{equation}
Note this estimation is independent from $ \varepsilon, m$. Thus,
from $\{v_{\varepsilon m}\}$ we can extract a subsequence
$\{v_{\varepsilon m_p}\}$ such that
\[
v_{\varepsilon m_p}\rightharpoonup v_\varepsilon\quad \mbox{weakly
in}\,\, X\,\, \mbox{as}\,\, p\rightarrow +\infty.
\]

By the assumptions of theorem
\begin{equation}\label{4.27}
a(v_{\varepsilon m_p},v_{\varepsilon m_p},u)+\langle
A_\varepsilon(v_{\varepsilon m_p}),u\rangle\rightarrow
a(v_\varepsilon,v_\varepsilon,u)+\langle
A_\varepsilon(v_\varepsilon) ,u\rangle\quad \forall u\in X.
\end{equation}
From \eqref{4.22}, \eqref{4.27}, \eqref{4.26} we know that
$v_\varepsilon$ is a solution to \eqref{4.21} and satisfies
\begin{equation}\label{4.28}
\|v_\varepsilon\|_X\leq \f{1}{\alpha}\|f\|_{X^*}.
\end{equation}
Subtracting the following two formula which are got from
\eqref{4.21}
\[
\begin{aligned}
& a(v_\varepsilon,v_\varepsilon,u)+\langle
A_\varepsilon(v_\varepsilon),u\rangle=\langle f, u\rangle\quad
\forall u\in X,\\&
a(v_\varepsilon,v_\varepsilon,v_\varepsilon)+\langle
A_\varepsilon(v_\varepsilon),v_\varepsilon\rangle=\langle f,
v_\varepsilon\rangle
\end{aligned}
\]
and taking into account that
\[
J_\varepsilon(u)-J_\varepsilon(v_\varepsilon)\geq \langle
A_\varepsilon(v_\varepsilon),u-v_\varepsilon\rangle
\]
which is due to convexity of $J_\varepsilon$, we come to the
following inequality
\begin{equation}\label{4.29}
a(v_\varepsilon,v_\varepsilon,u-v_\varepsilon)+J_\varepsilon(u)-J_\varepsilon(v_\varepsilon)\geq
\langle f, u-v_\varepsilon\rangle\quad \forall u\in X.
\end{equation}
By \eqref{4.28} we can choose $\{v_{\varepsilon_k}\}$ such that
\begin{equation}\label{4.30}
 v_{\varepsilon_k}\rightharpoonup
v^*\quad \mbox{weakly in}\,\,
X\,\,\mbox{as}\,\,\varepsilon_k\rightarrow 0.
\end{equation}
By virtue of assumption 1)
\[
|J_{\varepsilon_k}(v_{\varepsilon_k})-J(v^*)|\leq
|J_{\varepsilon_k}(v_{\varepsilon_k})-J(v_{\varepsilon_k})|+|J(v_{\varepsilon_k})-J(v^*)|\rightarrow
0 \quad\mbox{as}\,\,\varepsilon_k\rightarrow 0,
\]
and so
\begin{equation}\label{4.31}
J_{\varepsilon_k}(v_{\varepsilon_k})\rightarrow J(v^*)
\,\,\mbox{as}\,\,\varepsilon_k\rightarrow 0.
\end{equation}
Also
\begin{equation}\label{4.32}
J_{\varepsilon_k}(u)\rightarrow J(u)\quad \forall u\in X
\quad\mbox{as}\,\,\varepsilon_k\rightarrow 0.
\end{equation}
By virtue of assumption 2)
\begin{equation}\label{4.33}
\begin{aligned}
&a(v_{\varepsilon_k},v_{\varepsilon_k},u)\rightarrow a(v^*,v^*,u)
\quad \forall u\in X,\\& \liminf_{k\rightarrow \infty}
a(v_{\varepsilon_k},v_{\varepsilon_k},v_{\varepsilon_k})\geq
a(v^*,v^*,v^*).
\end{aligned}
\end{equation}
Taking into account \eqref{4.31}-\eqref{4.33}, from \eqref{4.29} we
get
\[
a(v^*,v^*,u-v^*)+J(u)-J(v^*)\geq \langle f, u-v^*\rangle\quad
\forall u\in X.
\]
By \eqref{4.28} we have
\begin{equation}\label{4.34}
 \|v^*\|_X\leq
\f{1}{\alpha}\|f\|_{X^*}.
\end{equation}
\hspace{12cm}$\square$\vspace*{.3cm}

\begin{remark}\label{r4.1}
The estimate of solutions in Theorem \ref{t4.1} is for all solutions
of the problem, but one in Theorem \ref{t4.3} is for the solution
guaranteed existence by the theorem.
\end{remark}

\section{Mixed boundary value problems of the Navier-Stokes and Stokes equations}
\setcounter{equation}{0}

In this section relying on the results in Section 4, we are
concerned with problems in Section 3.

\begin{theorem}\label{t5.1}
Let Assumption \ref{a3.1} hold, the surfaces
$\Gamma_{2j},\,\,\Gamma_{3j},\,\,\Gamma_{7j}$ be convex (cf.
Definition \ref{d2.1}),  $\alpha$ positive and
$\|U\|_{\mathbf{H}^1(\Omega)}$  small enough.  Then, when $f$ and
$\phi_i,\,i=2\sim7,$ are small enough, there exists a unique
solution to Problem I-VI for the stationary Navier-Stokes problem
with mixed boundary condition \eqref{2.8} in a neighborhood of $U$
in $\mathbf{H}^1(\Omega)$.
\end{theorem}
$\mathbf{Proof.}$\, Define a functional $J(u)$ by
 \eqref{3.7}, \eqref{3.8}. Trace operator is continuous
 and sum of convex functions is also convex.  Thus, the functional
  satisfies condition 1) of Theorem \ref{t4.2}.

Let $w=v-U$, $U$ be a function in Assumption \ref{a3.1} and
$a_{01}(\cdot, \cdot), a_{11}(\cdot,\cdot,\cdot)$ and $F_1\in
\textbf{V}(\Omega)^*$ be as \eqref{3.4}:
\begin{equation}\label{5.2}
\begin{aligned}
& a_{01}(w,u)=2\nu(\varepsilon(w),\varepsilon(u))+\langle(U\cdot
\nabla)w,u\rangle+\langle(w\cdot\nabla)U,u\rangle+2\nu(k(x)w,u)_{\Gamma_2}\\&
\hspace{1.6cm}+2\nu(S\tilde{w},\tilde{u})_{\Gamma_3}+2(\alpha(x)w,u)_{\Gamma_5}+\nu(k(x)w,u)_{\Gamma_7}\quad
\forall w,u\in \mathbf{V}(\Omega),
\\&a_{11}(w,u,v)=\langle (w\cdot\nabla)
u,v\rangle\quad \forall w,u,v\in \mathbf{V}(\Omega),\\
& \langle
F_1,u\rangle=-2\nu(\varepsilon(U),\varepsilon(u))-\langle(U\cdot\nabla)U,u\rangle-2\nu(k(x)U,u)_{\Gamma_2}-2\nu(S\tilde{U},\tilde{u})_{\Gamma_3}\\&
\hspace*{1.3cm}-2(\alpha(x)U,u)_{\Gamma_5}-\nu(k(x)U,u)_{\Gamma_7}
+\langle
f,u\rangle+\sum_{i=2,4,7}\langle\phi_i,u_n\rangle_{\Gamma_i}\\&\hspace*{1.3cm}+\sum_{i=3,5,6}\langle\phi_i,u\rangle_{\Gamma_i}\quad
\forall u\in \mathbf{V}(\Omega),
\end{aligned}
\end{equation}
By Korn's inequality
\begin{equation}\label{5.3}
2\nu(\varepsilon(w),\varepsilon(w))\geq \delta \|w\|^2_\mathbf{V}.
\end{equation}
On the other hand, applying H\"older inequality for $w\in
\mathbf{V}(\Omega)$ we have
\begin{equation}\label{5.4}
\begin{aligned}
|\langle(U\cdot
\nabla)w,w\rangle+\langle(w\cdot\nabla)U,w\rangle|\leq \gamma
\|w\|^2_{\mathbf{V}}\cdot\|U\|_{\mathbf{H}^1(\Omega)}.
\end{aligned}
\end{equation}

Therefore, if $\delta-\gamma\|U\|_{\mathbf{H}^1(\Omega)}=\beta_1>0$,
then by \eqref{5.3}, \eqref{5.4}, Assumption \ref{a3.1} and Lemma
\ref{l2.1} we have
\begin{equation}\label{5.5}
a_{01}(u,u)\geq \beta_1\|u\|^2_{\mathbf{V}}\quad \forall u\in
\textbf{V}(\Omega).
\end{equation}
It is easy to verify that
\begin{equation}\label{5.6}
|a_{01}(u,v)|\leq c\|u\|_{\mathbf{V}}\|v\|_{\mathbf{V}}\quad \forall
u,v\in \textbf{V}(\Omega).
\end{equation}
By \eqref {5.5} and \eqref{5.6}, $a_0(u,v)$ satisfies condition 2)
of Theorem \ref{t4.2}.

By H\"{o}lder inequality we can see
\begin{equation}\label{5.7}
|a_{11}(w,u,v)|\leq
c\|w\|_{\mathbf{V}}\|u\|_{\mathbf{V}}\|v\|_{\mathbf{V}}\quad \forall
w,u,v\in \textbf{V}(\Omega).
\end{equation}
which means $a_{11}(w,u,v)$ satisfies condition 3) of Theorem
\ref{t4.2}.

 Also
\begin{equation}\label{5.8}
\begin{aligned}
&\|F_1\|_{\mathbf{V}^*}\leq
M_1\Big(\|U\|_{\mathbf{H}^1}+\|U\|^2_{\mathbf{H}^1}+\|f\|_{\mathbf{V}^*}
+\sum_{i=2,4,7}\|\phi_i\|_{H^{-\f{1}{2}}(\Gamma_i)}+\sum_{i=3,5,6}\|\phi_i\|_{\mathbf{H}^{-\f{1}{2}}(\Gamma_i)}\Big),
\end{aligned}
\end{equation}
where $M_1$ depends on  mean curvature of $\Gamma_7$, shape operator
of $\Gamma_3$, $\nu$ and $\alpha$.

By Theorem \ref{t4.2}, if $\|U\|_{\mathbf{H}^1},
\|f\|_{\mathbf{V}^*}, \|\phi_i\|_{H^{-\f{1}{2}}(\Gamma_i)},
i=2,4,7,$ and $\|\phi_i\|_{\mathbf{H}^{-\f{1}{2}}(\Gamma_i)},
i=3,5,6,$ are small enough, then there exists a unique solution
$w\in K(\Omega)$ to
\begin{equation}\label{5.9}
\begin{aligned}
& a_{01}(w,u-w)+a_{11}(w,w,u-w)+J(u)-J(w)\geq\langle
F_1,u-w\rangle\quad \forall u\in K(\Omega).
\end{aligned}
\end{equation}
Since $v=w+U$ is solution, we come to the asserted conclusion.
$\square$\vspace*{.3cm}

\begin{theorem}\label{t5.2}
Let Assumption \ref{a3.1} hold, the surfaces
$\Gamma_{2j},\,\,\Gamma_{3j},\,\,\Gamma_{7j}$ be convex, $\alpha$
positive and $\|U\|_{\mathbf{H}^1(\Omega)}$ small enough.  Then, for
any $f$
 $\phi_i,\,i=2\sim7,$ there exists a solution $v$ to Problem II-VI for the stationary Navier-Stokes problem
with mixed boundary condition \eqref{2.9} in a neighborhood of $U$
in $\mathbf{H}^1(\Omega)$ and all solutions satisfy
\begin{equation}\label{5.10}
\begin{aligned}
\|v-U\|_{\mathbf{H}^1}\leq\f{
M_1}{\delta-\gamma\|U\|_{\mathbf{H}^1}}\Big(&\|U\|_{\mathbf{H}^1}+\|U\|^2_{\mathbf{H}^1}+\|f\|_{\mathbf{V}^*}\\&
+\sum_{i=2,4,7}\|\phi_i\|_{H^{-\f{1}{2}}(\Gamma_i)}+\sum_{i=3,5,6}\|\phi_i\|_{\mathbf{H}^{-\f{1}{2}}(\Gamma_i)}\Big),
\end{aligned}
\end{equation}
where $\delta, \gamma, M_1$ are as \eqref{5.13}, \eqref{5.14},
\eqref{5.22}.

If $\|U\|_{\mathbf{H}^1}, \|f\|_{\mathbf{V}^*},
 \|\phi_i\|_{H^{-\f{1}{2}}(\Gamma_i)}, i=2,4,7,$ and
$\|\phi_i\|_{\mathbf{H}^{-\f{1}{2}}(\Gamma_i)}, i=3,5,6,$ are small
enough, then the solution is unique.
\end{theorem}
$\mathbf{Proof.}$\, Define a functional $J(u)$ by
 \eqref{3.7}, \eqref{3.8}. Then, this functional
satisfies condition 1) of Theorem \ref{t4.1}.

Let $a_{02}(\cdot, \cdot), a_{12}(\cdot,\cdot,\cdot)$ and $F_2\in
V^*$ are as \eqref{3.18}:
\begin{equation}\label{5.12}
\begin{aligned}
&a_{02}(w,u)=2\nu(\varepsilon(w),\varepsilon(u))+\langle\mbox{rot}\,U\times
w,u\rangle+\langle\mbox{rot}\,w\times
U,u\rangle+2\nu(k(x)w,u)_{\Gamma_2}\\&
\hspace{2.1cm}+2\nu(S\tilde{w},\tilde{u})_{\Gamma_3}+2(\alpha(x)w,u)_{\Gamma_5}+\nu(k(x)w,u)_{\Gamma_7},
\\&a_{12}(w,u,v)=\langle\mbox{rot}\,w\times
u,v\rangle,\\& \langle
F_2,u\rangle=-2\nu(\varepsilon(U),\varepsilon(u))-\langle\mbox{rot}\,U\times
U,u\rangle-2\nu(k(x)U,u)_{\Gamma_2}-2\nu(S\tilde{U},\tilde{u})_{\Gamma_3}\\
&\hspace*{1cm}-2(\alpha(x)U,u)_{\Gamma_5}-\nu(k(x)U,u)_{\Gamma_7}+\langle
f,u\rangle+\sum_{i=2,4,7}\langle\phi_i,u_n\rangle_{\Gamma_i}+\sum_{i=3,5,6}\langle\phi_i,u\rangle_{\Gamma_i}.
\end{aligned}
\end{equation}
By Korn's inequality
\begin{equation}\label{5.13}
2\nu(\varepsilon(w),\varepsilon(w))\geq \delta \|w\|^2_\mathbf{V}.
\end{equation}
On the other hand, for any $w\in \mathbf{V}(\Omega)$ we have
\begin{equation}\label{5.14}
\begin{aligned}
&\langle\mbox{rot}\,U\times w,w\rangle=0,\\&
|\langle\mbox{rot}\,w\times U,w\rangle|\leq \gamma
\|w\|^2_{\mathbf{V}}\cdot\|U\|_{\mathbf{H}^1(\Omega)}.
\end{aligned}
\end{equation}
Therefore, if $\delta-\gamma\|U\|_{\mathbf{H}^1(\Omega)}=\beta_1>0$,
then by \eqref{5.13}, \eqref{5.14}, Assumption \ref{a3.1} and Lemma
\ref{l2.1} we have
\begin{equation}\label{5.15}
a_{02}(u,u)\geq \beta_1\|u\|^2_{\mathbf{V}}\quad \forall u\in
\textbf{V}(\Omega).
\end{equation}
It is easy to verify
\begin{equation}\label{5.16}
|a_{02}(u,v)|\leq
c\|u\|_{\mathbf{V}(\Omega)}\|v\|_{\mathbf{V}(\Omega)}\quad \forall
u,v\in \textbf{V}(\Omega).
\end{equation}
Then, \eqref{5.15}, \eqref{5.16} show that $a_{02}(u,v)$ satisfy
condition 2) of Theorem \ref{t4.1}.

 By a property of mixed product,
\begin{equation}\label{5.17}
a_{12}(w,u,u)=\langle\mbox{rot}\,w\times u,u\rangle=0\quad \forall
w\in\textbf{V}^{\f{2}{3}}(\Omega), \forall u\in\mathbf{V}(\Omega),
\end{equation}
where $\textbf{V}^{\f{2}{3}}(\Omega)=\{u\in
\mathbf{H}^{\f{2}{3}}(\Omega):\mbox{div}\,u=0,\, u|_{\Gamma_1}=0,\,
u_\tau|_{(\Gamma_2\cup\Gamma_4\cup\Gamma_7\cup\Gamma_9)}=0,\, u\cdot
n|_{(\Gamma_3\cup\Gamma_5\cup\Gamma_8)}=0\}$.
 On the other hand, by density argument we get
\begin{equation}\label{5.18}
a_{12}(w,u,v)=\langle\mbox{rot}\,w\times
u,v\rangle=-\langle\mbox{rot}\,w, v\times u\rangle.
\end{equation}
When $u,v\in \textbf{V}(\Omega)$, $v\times u\in
\textbf{H}^{\f{1}{2}}(\Omega)$ and
\begin{equation}\label{5.19}
\|v\times u\|_{ \textbf{H}^{\f{1}{3}}(\Omega)}\leq c_1\|v\times
u\|_{ \textbf{H}^{\f{1}{2}}(\Omega)}\leq c\|v\|_{
\textbf{V}(\Omega)}\|u\|_{\textbf{V}(\Omega)}.
\end{equation}
(cf. Theorem 1.4.4.2 in \cite{g1}.) Also, if
$w\in\textbf{V}^{\f{2}{3}}(\Omega)$, then
$\mbox{rot}\,w\in\textbf{H}^{-\f{1}{3}}(\Omega)$ and
\begin{equation}\label{5.20}
\|\mbox{rot}\,w\|_{ \textbf{H}^{-\f{1}{3}}(\Omega)}\leq c\|w\|_{
\textbf{H}^{\f{2}{3}}(\Omega)}.
\end{equation}
(cf. Proposition 12.1, ch. 1 in \cite{lm}.) Since
$H^{\f{1}{3}}_0(\Omega)=H^{\f{1}{3}}(\Omega)$(cf. Theorem 11.1, ch.
1 in \cite{lm}), by \eqref{5.18}-\eqref{5.20} we get
\begin{equation}\label{5.21}
|a_{12}(w,u,v)|\leq
K\|w\|_{\textbf{V}^{\f{2}{3}}(\Omega)}\|u\|_{\textbf{V}(\Omega)}\|v\|_{\textbf{V}(\Omega)}\quad
\forall w\in\textbf{V}^{\f{2}{3}}(\Omega), \forall
u,v\in\mathbf{V}(\Omega).
\end{equation}
Since
$\textbf{V}(\Omega)\hookrightarrow\hookrightarrow\textbf{V}^{\f{2}{3}}(\Omega)$,
setting $X=\textbf{V}(\Omega),$ $X_1=\textbf{V}^{\f{2}{3}}(\Omega)$
by \eqref{5.17}, \eqref{5.21} $a_{11}(w,u,v)$ satisfies condition 3)
of Theorem \ref{t4.1}.

 Also, we have
\begin{equation}\label{5.22}
\begin{aligned}
&\|F_2\|_{\mathbf{V}^*}\leq
M_1\Big(\|U\|_{\mathbf{H}^1}+\|U\|^2_{\mathbf{H}^1}+\|f\|_{\mathbf{V}^*}
+\sum_{i=2,4,7}\|\phi_i\|_{H^{-\f{1}{2}}(\Gamma_i)}+\sum_{i=3,5,6}\|\phi_i\|_{\mathbf{H}^{-\f{1}{2}}(\Gamma_i)}\Big),
\end{aligned}
\end{equation}
where $M_1$ depends on  mean curvature, shape operator, $\nu$ and
$\alpha$.

Therefore, by Theorem \ref{t4.1}, we have existence and  an estimate
of solutions to
\[
\begin{aligned}
& a_{02}(w,u-w)+a_{12}(w,w,u-w)+J(u)-J(w)\geq\langle
F_2,u-w\rangle\quad \forall u\in \mathbf{V}(\Omega).
\end{aligned}
\]
Since $v=w+U$ is solution to the given problem, we have existence of
solutions and the estimate  \eqref{5.10}.

If $\|U\|_{\mathbf{H}^1}, \|f\|_{\mathbf{V}^*},$ $
\|\phi_i\|_{H^{-\f{1}{2}}(\Gamma_i)}, i=2,4,7,$ and
$\|\phi_i\|_{\mathbf{H}^{-\f{1}{2}}(\Gamma_i)}, i=3,5,6,$ are small
enough, then the solution is unique. $\square$\vspace*{.3cm}

Let us consider a special case of the Navier-Stokes problem with
boundary condition \eqref{2.8} in which there is not any flux across
boundary except $\Gamma_1, \Gamma_8$.

\begin{theorem}\label{t5.2.0}
Let Assumption \ref{a3.1} hold, $\Gamma_i=\emptyset
(i=2,4,6,7,9-11),$ the surfaces $\Gamma_{3j}$ be convex, $\alpha$
positive and $\|U\|_{\mathbf{H}^1(\Omega)}$ small enough. Then, for
any $f$ and $\phi_i,\,i=3,5$ there exists a solution $v$ to Problem
I-VI for the stationary Navier-Stokes problem with mixed boundary
condition \eqref{2.8} and all solutions satisfy
\begin{equation}\label{5.23}
\begin{aligned}
&\|v-U\|_{\mathbf{H}^1}\leq\f{
M_1}{\delta-\gamma\|U\|_{\mathbf{H}^1}}\Big(\|U\|^2_{\mathbf{H}^1}+\|f\|_{\mathbf{V}^*}
+\sum_{i=3,5}\|\phi_i\|_{\mathbf{H}^{-\f{1}{2}}(\Gamma_i)}\Big),
\end{aligned}
\end{equation}
where $\delta, \gamma, M_1$ are as \eqref{5.3}, \eqref{5.4},
\eqref{5.8}.

In addition, if $\|f\|_{\mathbf{V}^*},$ $
\|\phi_i\|_{\mathbf{H}^{-\f{1}{2}}(\Gamma_i)}, i=3,5,$ are small
enough, then the solution is unique.
\end{theorem}
$\mathbf{Proof.}$\, Define a functional $J(u)=j_\tau(u)$ by
 \eqref{3.7}, \eqref{3.8}. Then, the functional
  satisfies condition 1) of Theorem \ref{t4.2}.

Let $w=v-U$, $U$ be a function in Assumption \ref{a3.1} and
$a_{01}(\cdot, \cdot), a_{11}(\cdot,\cdot,\cdot)$ and $F_1\in
\textbf{V}(\Omega)^*$ be as \eqref{3.4}:
\[
\begin{aligned}
& a_{01}(w,u)=2\nu(\varepsilon(w),\varepsilon(u))+\langle(U\cdot
\nabla)w,u\rangle+\langle(w\cdot\nabla)U,u\rangle+2\nu(S\tilde{w},\tilde{u})_{\Gamma_3}+2(\alpha(x)w,u)_{\Gamma_5}\\&\hspace{9cm}\quad
\forall w,u\in \mathbf{V}(\Omega),
\\&a_{11}(w,u,v)=\langle (w\cdot\nabla)
u,v\rangle\quad \forall w,u,v\in \mathbf{V}(\Omega),\\
& \langle
F_1,u\rangle=-2\nu(\varepsilon(U),\varepsilon(u))-\langle(U\cdot\nabla)U,u\rangle-2\nu(S\tilde{U},\tilde{u})_{\Gamma_3}-2(\alpha(x)U,u)_{\Gamma_5}
+\langle
f,u\rangle+\sum_{i=3,5}\langle\phi_i,u\rangle_{\Gamma_i}\\&\hspace*{9cm}\quad
\forall u\in \mathbf{V}(\Omega),
\end{aligned}
\]
We can see that the condition 2) in Theorem \ref{t4.1} are
satisfied(cf. proof of Theorem \ref{t5.1}).

 By the condition of theorem,
\begin{equation}\label{5.22.0}
a_{11}(w,u,u)=\langle (w\cdot\nabla) u,u\rangle=0\quad \forall
w\in\textbf{V}^{\f{2}{3}}(\Omega), \forall u\in \mathbf{V}(\Omega).
\end{equation}
By H\"{o}lder inequality we can see
\begin{equation}\label{5.23.0}
|a_{11}(w,u,v)|\leq
K\|w\|_{\textbf{V}^{\f{2}{3}}(\Omega)}\|u\|_{\textbf{V}(\Omega)}\|v\|_{\textbf{V}(\Omega)}\quad
\forall w\in\textbf{V}^{\f{2}{3}}(\Omega), \forall
u,v\in\mathbf{V}(\Omega).
\end{equation}
By \eqref{5.22.0}, \eqref{5.23.0}, $a_{11}(w,u,v)$ satisfies
condition 3) of Theorem \ref{t4.1}.

Applying Theorem \ref{t4.1} to
\[
\begin{aligned}
& a_{01}(w,u-w)+a_{11}(w,w,u-w)+J(u)-J(w)\geq\langle
F_1,u-w\rangle\quad \forall u\in K(\Omega),
\end{aligned}
\]
we come to the asserted conclusion. $\square$\vspace*{.3cm}

\begin{remark}\label{r5.1}
Assumption $\Gamma_i=\varnothing, i=2,4,6,7,9-11,$ is only used
to get \eqref{5.22.0}.
\end{remark}\vspace*{.3cm}

Relying on Theorem \ref{t4.3}, again let us study the problem
concerned in Theorem \ref{t5.2.0}. This is generalization of methods
used in previous papers relying on smooth approximation of
functional in variational inequalities(cf. \cite{ll3}).

\begin{lemma}\label{l5.1} Let $X, Y$ be reflex Banach spaces, an operator $i\in (X\rightarrow Y)$
be completely linear continuous, $j\in(Y\rightarrow R)$ be convex
and Gateaux derivative $Dj(y)= a(y)$ for $y\in Y$. Then, $J(v)\equiv
j(iv)\in (X\rightarrow R)$ is convex, $DJ(v)\equiv A(v)=i^*a(iv)$,
where $i^*$ is the operator adjoint to $i$, and $A\in (X\rightarrow
X^*)$ is weak continuous.
\end{lemma}
$\mathbf{Proof.}$\,It is easy to verify convexity of $ J$.
 \[
\begin{aligned}
\langle A(v), u\rangle_X&=\lim_{t\rightarrow
0}\f{J(v+tu)-J(u)}{t}=\lim_{t\rightarrow
0}\f{j(i(v+tu))-j(iu)}{t}\\&=\langle a(iv), iu\rangle_Y=\langle
i^*a(iv), u\rangle_X\qquad \forall v,u\in X,
\end{aligned}
\]
which means $A(v)=i^*a(iv)$.

Let $v_n\rightharpoonup  v$ weakly in $X$. Since Gateaux derivative
of a finite convex functional is monotone and demi-continuous(cf.
Lemmas 4.10, 4.12, ch. 3 in \cite{ggz}) and $iv_n\rightarrow iv$ in
$Y$,
\[
\begin{aligned}
\langle A(v_n), u\rangle_X&=\langle i^*a(iv_n), u\rangle_X=\langle
a(iv_n), iu\rangle_Y\rightarrow\langle a(iv), iu\rangle_Y=\langle
i^*a(iv), u\rangle_X\qquad \forall u\in X,
\end{aligned}
\]
that is, $DJ=A\in(X\rightarrow X^*)$ is weak continuous.
$\square$\vspace*{.3cm}

\begin{theorem}\label{t5.3}
Let Assumption \ref{a3.1} hold, $\Gamma_i=\emptyset
(i=2,4,6,7,9-11),$ the surfaces $\Gamma_{3j}$ be convex, $\alpha$
positive and $\|U\|_{\mathbf{H}^1(\Omega)}$ small enough. Then, for
any $f$ and $\phi_i,\,i=3,5,$ there exists a solution $v$ to Problem
I-VI for the stationary Navier-Stokes problem with mixed boundary
condition \eqref{2.8} and the solution satisfies the estimate \eqref{5.23}.
\end{theorem}
$\mathbf{Proof.}$\,
 Define an operator $i\in \left(\mathbf{V}(\Omega)\rightarrow
 \textbf{L}_\tau^2(\Gamma_8)\right)$ by $iu=u|_{\Gamma_8}$ and a functional $J\in \left(\mathbf{V}(\Omega)\rightarrow
R\right)$ by $J(v)\equiv j_\tau(iv)$, where $j_\tau$ is as
\eqref{3.7}. Since the trace operator
$(\textbf{V}(\Omega)\rightarrow H^{\f{1}{2}}(\partial \Omega))$ is
continuous and $H^{\f{1}{2}}(\partial
\Omega)\hookrightarrow\hookrightarrow L^2(\partial \Omega)$, the
operator $i$ is compact, and by Lemma \ref{l5.1}
$J\in\left(\mathbf{V}(\Omega)\rightarrow R\right) $ is weak
continuous and convex.

Define a functional $J_\varepsilon\in
\left(\mathbf{V}(\Omega)\rightarrow R\right)$ by
\begin{equation}\label{5.24}
\begin{aligned}
&J_\varepsilon(v)= j_{\tau\varepsilon}(iv),\\ &
j_{\tau\varepsilon}(\eta)=\int_{\Gamma_8}g_\tau\rho_\varepsilon(\eta)\,ds,\\&
\rho_\varepsilon(\eta)=\begin{cases}|\eta|-\varepsilon/2&|\eta|>\varepsilon,\\
|\eta|^2/2\varepsilon&|\eta|\leq\varepsilon\,.\end{cases}
\end{aligned}
\end{equation}
Since
\[
|j_{\tau\varepsilon}(\eta)-j_\tau(\eta)|\leq
\f{\varepsilon}{2}|g_\tau| \quad \forall
\eta\in\textbf{L}^2_\tau(\Gamma_8)
\]
(cf. Lemma 2.1 in \cite{ll3}), we have
\begin{equation}\label{5.25}
|J_\varepsilon(v)-J(v)|\leq \f{\varepsilon}{2}|g_\tau| \quad \forall
v\in\textbf{V}(\Omega).
\end{equation}
Also, $j_{\tau\varepsilon}$ is convex, and so its Gateaux derivative
is demi-continuous. Thus, by Lemma \ref{l5.1} $DJ_\varepsilon\equiv
A_\varepsilon\in(\mathbf{V}(\Omega)\rightarrow\mathbf{V}(\Omega)^*)$
is weak continuous. By this fact together \eqref{5.25}, condition 1)
of Theorem \ref{t4.3} is satisfied.

Under Assumption of theorem $a_{01}(\cdot, \cdot),
a_{11}(\cdot,\cdot,\cdot)$ and $F_1\in V^*$ of \eqref{3.4} are as
follows.
\begin{equation}\label{5.26}
\begin{aligned}
& a_{01}(u,v)=2\nu(\varepsilon(u),\varepsilon(v))+\langle(U\cdot
\nabla)u,v\rangle+\langle(u\cdot\nabla)U,v\rangle\\&
\hspace{1.6cm}+2\nu(S\tilde{u},\tilde{v})_{\Gamma_3}+2(\alpha(x)u,v)_{\Gamma_5}\quad
\forall u,v\in \mathbf{V}(\Omega),
\\&a_{11}(w,u,v)=\langle (w\cdot\nabla)
u,v\rangle\quad \forall w,u,v\in \mathbf{V}(\Omega),\\
& \langle
F_1,u\rangle=-2\nu(\varepsilon(U),\varepsilon(u))-\langle(U\cdot\nabla)U,u\rangle-2\nu(S\tilde{U},\tilde{u})_{\Gamma_3}\\&
\hspace*{1.3cm}-2(\alpha(x)U,u)_{\Gamma_5} +\langle
f,u\rangle+\sum_{i=3,5}\langle\phi_i,u\rangle_{\Gamma_i}\quad
\forall u\in \mathbf{V}(\Omega),
\end{aligned}
\end{equation}
By Korn's inequality
\begin{equation}\label{5.27}
2\nu(\varepsilon(u),\varepsilon(u))\geq \delta \|u\|^2_\mathbf{V}.
\end{equation}

On the other hand, for any $w\in \mathbf{V}(\Omega)$ we have
\begin{equation}\label{5.28}
\begin{aligned}
|\langle(U\cdot
\nabla)u,u\rangle+\langle(u\cdot\nabla)U,u\rangle|\leq \gamma
\|u\|^2_{\mathbf{V}}\cdot\|U\|_{\mathbf{H}^1(\Omega)}.
\end{aligned}
\end{equation}

Therefore, if $\delta-\gamma\|U\|_{\mathbf{H}^1(\Omega)}=\beta_1>0$,
then by \eqref{5.27}, \eqref{5.28}, Assumption \ref{a3.1} and Lemma
\ref{l2.1} we have
\begin{equation}\label{5.29}
a_{01}(u,u)\geq \beta_1\|u\|^2_{\mathbf{V}}\quad \forall u\in
\textbf{V}(\Omega).
\end{equation}
Under condition $\Gamma_i=\emptyset, i=2,4,6,7,9,10,11,$ it is easy
to verify that
\begin{equation}\label{5.30}
a_{11}(v,v,v)=0\quad \forall v\in \textbf{V}(\Omega).
\end{equation}
Let
\[
a(w,u,v)=a_{01}(u,v)+a_{11}(w,u,v).
\]
Then, by \eqref{5.29}, \eqref{5.30} we have
\begin{equation}\label{5.31}
a(v,v,v)\geq \beta_1\|u\|^2_{\mathbf{V}}\quad \forall v\in
\textbf{V}(\Omega).
\end{equation}

Let us prove that when $v_m \rightharpoonup v$ weakly in
$\textbf{V}(\Omega)$, for a subsequence $\{v_{m_p}\}$
\begin{equation}\label{5.32}
a(v_{m_p},v_{m_p},u)\rightarrow a(v;v,u) \quad\forall
u\in\textbf{V}(\Omega).
\end{equation}
To this end, first let us prove that when $v_m \rightharpoonup v$
weakly in $\textbf{V}(\Omega)$, for a subsequence $\{v_{m_p}\}$
\begin{equation}\label{5.33}
 a_{01}(v_m,u)\rightarrow  a_{01}(v_m,u)\quad \forall u\in
\textbf{V}(\Omega).
\end{equation}
Since $U_iu_j\in L^2(\Omega), i,j=1-3$, and $\p_i v_m\rightharpoonup
v$ in $\text{L}(\Omega)^2$, we have
\begin{equation}\label{5.34}
\langle(U\cdot \nabla)v_m,u\rangle\rightarrow\langle(U\cdot
\nabla)v,u\rangle\quad \text{as}\,\,m\rightarrow\infty.
\end{equation}
By H\"older inequalities
 $$
 |\langle((v_m-v)\cdot\nabla)U,u\rangle|\leq c \|v_m-v\|_{\textbf{L}^3(\Omega)}\|\nabla
 U\|_{\textbf{L}^2(\Omega)}\|u\|_{\textbf{L}^6(\Omega)}.
 $$
Since $H^1(\Omega)\hookrightarrow\hookrightarrow L^3(\Omega)$, we
can choose a subsequence $\{v_{m_p}\}$ such that
 $v_{mp}\rightarrow v$ in $\textbf{L}^3(\Omega)$. Then, we have
\begin{equation}\label{5.35}
\langle(v_{m_p}\cdot\nabla)U,u\rangle\rightarrow
\langle(v\cdot\nabla)U,u\rangle\quad\,\,\text{as}\,\,
m_p\rightarrow\infty.
\end{equation}
It is easy to verify convergence of other terms. Thus, using
\eqref{5.34}, \eqref{5.35}, we have \eqref{5.33}.

Using H\"older inequality and $a_{11}(v,u,w)=-a_{11}(v,w,u)$, we
have
\[
\begin{aligned}
&|a_{11}(v_m,v_m,u)-a_{11}(v,v,u)|\\&\leq|a_{11}(v_m,v_m,u)-a_{11}(v,v_m,u)|+|a_{11}(v,v_m,u)-a_{11}(v,v,u)|\quad
\\&\leq c\big(\|v_m-v\|_{\textbf{L}^3(\Omega)}\|\nabla v_m\|_{\textbf{L}^2(\Omega)}\|u\|_{\textbf{L}^6(\Omega)}+\|v\|_{\textbf{L}^6(\Omega)}\|\nabla u\|_{\textbf{L}^2(\Omega)}\|v_m-v\|_{\textbf{L}^3(\Omega)}\big)\quad\forall u\in \mathbf{V}(\Omega).
\end{aligned}
\]
Thus,  we have
\begin{equation}\label{5.36}
a_{11}(v_{m_p},v_{m_p},u)\rightarrow a_{11}(v,v,u)\quad\forall u\in
\mathbf{V}(\Omega)\quad\text{as}\,\, m_p\rightarrow\infty .
\end{equation}
From \eqref{5.33}, \eqref{5.36} we get \eqref{5.32}.

Let us prove that
\begin{equation}\label{5.37}
\liminf_{m\rightarrow\infty} a(v_{m_p},v_{m_p},v_{m_p})\geq
a(v,v,v).
\end{equation}
By lower semi-continuity of norm
\begin{equation}\label{5.38}
\liminf_{m\rightarrow\infty}2\nu(\varepsilon(v_m),\varepsilon(v_m))\geq2\nu(\varepsilon(v),\varepsilon(v))
\quad \text{as}\,\,v_m\rightharpoonup
v\,\,\text{in}\,\,\textbf{V}(\Omega).
\end{equation}
It is easy to prove that
\begin{equation}\label{5.39}
2\nu(S\tilde{v}_{m},\tilde{u})_{\Gamma_3}+2(\alpha(x)v_{m},u)_{\Gamma_5}
\rightarrow
2\nu(S\tilde{v},\tilde{u})_{\Gamma_3}+2(\alpha(x)v,u)_{\Gamma_5}
\quad \forall u\in \textbf{V}(\Omega).
\end{equation}
Using H\"older inequality and $a_{11}(v,v_m,u)=-a_{11}(v,u,v_m)$, we
have
\[
\begin{aligned}
|a_{11}&(v_m,v_m,v_m)-a_{11}(v,v,v)|\\&\leq|a_{11}(v_m,v_m,v_m)-a_{11}(v,v_m,v_m)|+|a_{11}(v,v_m,v_m)-a_{11}(v,v_m,v)|\\
&\hspace{1cm}+|a_{11}(v,v_m,v)-a_{11}(v,v,v)|\quad
\\&\leq c\big(\|v_m-v\|_{\textbf{L}^3(\Omega)}\|\nabla v_m\|_{\textbf{L}^2(\Omega)}\|v_m\|_{\textbf{L}^6(\Omega)}+\|v\|_{\textbf{L}^6(\Omega)}\|\nabla v_m\|_{\textbf{L}^2(\Omega)}\|v_m-v\|_{\textbf{L}^3(\Omega)}\\
&\hspace{1cm}+\|v\|_{\textbf{L}^6(\Omega)}\|\nabla
v\|_{\textbf{L}^2(\Omega)}\|v_m-v\|_{\textbf{L}^3(\Omega)}\big),
\end{aligned}
\]
which implies
\begin{equation}\label{5.39.1}
a_{11}(v_{m_p},v_{m_p},v_{m_p})\rightarrow
a_{11}(v,v,v)\quad\text{as}\,\, m_p\rightarrow\infty.
\end{equation}
From \eqref{5.37}-\eqref{5.39.1}, we have \eqref{5.37}.

By virtue of \eqref{5.31}, \eqref{5.32} and \eqref{5.37}, condition
2) of Theorem \ref{t4.3} is satisfied. Therefore, by Theorem
\ref{t4.3} we have existence of a solution $w\in \mathbf{V}(\Omega)$
to
\begin{equation}\label{5.40}
\begin{aligned}
& a_{01}(w,u-w)+a_{11}(w,w,u-w)+j_\tau(u)-j_\tau(w)\geq\langle
F_1,u-w\rangle\quad \forall u\in \mathbf{V}(\Omega)
\end{aligned}
\end{equation}
and an estimate. Since $v=w+U$ is a solution, we come to the
asserted conclusion.
 $\square$\vspace*{.3cm}

\begin{remark}\label{r5.2}
The estimate of solution of Theorem \ref{t5.3} is not for all
solutions, and so Theorem \ref{t5.3} is weaker than Theorem
\ref{t5.2.0}.
\end{remark}

Let us consider Problem III for the Stokes system.

\begin{theorem}\label{t5.4}
Let Assumption \ref{a3.1} hold, the surfaces
$\Gamma_{2j},\,\,\Gamma_{3j},\,\,\Gamma_{7j}$ be convex and $\alpha$
positive. Then, for any $f$
 $\phi_i,\,i=2\sim7,$ there exists a unique solution $v$ to Problem III-VI for the stationary Stokes problem
with mixed boundary condition \eqref{2.8} and
\begin{equation}\label{5.41}
\begin{aligned}
&\|v-U\|_{\mathbf{H}^1}\leq\f{
M_1}{\delta}\Big(\|U\|_{\mathbf{H}^1}+\|f\|_{\mathbf{V}^*}
+\sum_{i=2,4,7}\|\phi_i\|_{H^{-\f{1}{2}}(\Gamma_i)}+\sum_{i=3,5,6}\|\phi_i\|_{\mathbf{H}^{-\f{1}{2}}(\Gamma_i)}\Big),
\end{aligned}
\end{equation}
where $\delta,  M_1$ are as \eqref{5.13}, \eqref{5.22}(for $F_3$
instead of $F_2$).

If $v_1, v_2$ are solutions, respectively, to Problem-III-VI with
$g_{\tau1}, g_{n1}, g_{+n1}, g_{-n1}, f_1, h_i^1, \phi_i^1$ and
$g_{\tau2}, g_{n2}, g_{+n2}, g_{-n2}, $ $ f_2, h_i^1, \phi_i^2$,
then
\begin{equation}\label{5.42}
\begin{aligned}
\|v_1-v_2\|_{\mathbf{H}^1}\leq\f{
M_1}{\delta}\Big(&\|U_1-U_2\|_{\mathbf{H}^1}+\|f_1-f_2\|_{\mathbf{V}^*}
+\|g_{\tau1}-g_{\tau2}\|_{\texttt{L}_\tau^2(\Gamma_8)}\\&+\|g_{n1}-g_{n2}\|_{L^2(\Gamma_9)}+\|g_{+n1}-g_{+n2}\|_{L^2(\Gamma_{10})}+\|g_{-n1}-g_{-n2}\|_{L^2(\Gamma_{10})}\\&+\sum_{i=2,4,7}\|\phi_i^1-\phi_i^2\|_{H^{-\f{1}{2}}(\Gamma_i)}+\sum_{i=3,5,6}\|\phi_i^1-\phi_i^2\|_{\mathbf{H}^{-\f{1}{2}}(\Gamma_i)}\Big)+\|U_1-U_2\|_{\mathbf{H}^1},
\end{aligned}
\end{equation}
where $U_j, j=1,2,$ are the functions in Assumption \ref{a3.1} with
$h_i^j$ instead $h_i.$
\end{theorem}

 $\mathbf{Proof.}$\, By arguments similar to proof of
Theorem \ref{t4.2} we can apply the well known result for
variational inequality
\begin{equation}\label{5.43} a_{03}(w,u-w)+J(u)-J(w)\geq
\langle F_3, u-w\rangle\quad \forall u\in X,
\end{equation}
where $J(u)$ is defined by \eqref{3.7}, \eqref{3.8} and
$a_{03}(v,u), F_3$ are as \eqref{3.24}:
\[
\begin{aligned}
&
a_{03}(w,u)=2\nu(\varepsilon(w),\varepsilon(u))+2\nu(k(x)w,u)_{\Gamma_2}\\&
\hspace{1.6cm}+2\nu(S\tilde{w},\tilde{u})_{\Gamma_3}+2(\alpha(x)w,u)_{\Gamma_5}+\nu(k(x)w,u)_{\Gamma_7}\quad
\forall w,u\in \mathbf{V}(\Omega),\\
& \langle
F_3,u\rangle=-2\nu(\varepsilon(U),\varepsilon(u))-2\nu(k(x)U,u)_{\Gamma_2}-2\nu(S\tilde{U},\tilde{u})_{\Gamma_3}
-2(\alpha(x)U,u)_{\Gamma_5}\\& \hspace*{1.3cm}
-\nu(k(x)U,u)_{\Gamma_7}+\langle
f,u\rangle+\sum_{i=2,4,7}\langle\phi_i,u_n\rangle_{\Gamma_i}+\sum_{i=3,5,6}\langle\phi_i,u\rangle_{\Gamma_i}\quad
\forall u\in \mathbf{V}(\Omega).
\end{aligned}
\]
Thus, we have a unique existence of solution  and estimate
\eqref{5.41}.

If $v_1=w_1+U_1, v_2=w_2+U_2$ are solutions corresponding to the
given data, we get
\begin{equation}\label{5.44}
\begin{aligned}
&a_{03}(w_1,u-w_1)+J_1(u)-J_1(w_1)\geq \langle F_3^1,
u-w_1\rangle,\\
& a_{03}(w_2,u-w_2)+J_2(u)-J_2(w_2)\geq \langle F_3^2,
u-w_2\rangle\quad \forall u\in \textbf{V}(\Omega),
\end{aligned}
\end{equation}
where $J_j(u), F_3^j, j=1,2,$ are one corresponding to $U_j, g_{\tau
j}, g_{nj},  g_{+nj}, g_{-nj},f_j, h_i^j, \phi_i^j$. Putting $u=w_2,
u=w_1,$ respectively, in the first and second one in \eqref{5.44}
and adding those, we have
\begin{equation}\label{5.45}
\begin{aligned}
&a_{03}(w_1-w_2,w_2-w_1)+J_1(w_2)-J_1(w_1)+J_2(w_1)-J_2(w_2)\geq
\langle F_3^1-F_3^2, w_2-w_1\rangle.
\end{aligned}
\end{equation}
By Korn's inequality and Lemma \ref{l2.1} we have
\begin{equation}\label{5.46}
a_{03}(w_1-w_2,w_1-w_2)\geq \delta \|w_1-w_2\|^2_\mathbf{V}.
\end{equation}
From \eqref{5.45}, \eqref{5.46} we have
\begin{equation}\label{5.47}
\begin{aligned}
\|w_1-w_2\|^2_{\mathbf{V}}\leq\f{ 1}{\delta}\Big(|\langle
F_3^1-F_3^2,
w_2-w_1\rangle|+|J_1(w_2)-J_1(w_1)+J_2(w_1)-J_2(w_2)|\Big).
\end{aligned}
\end{equation}
Since $w_1, w_2\in K(\Omega)$,
\begin{equation}\label{5.48}
\begin{aligned}
&J_1(w_2)-J_1(w_1)=\int_{\Gamma_8}g_{\tau
1}(|w_{2\tau}|-|w_{1\tau}|)\,ds+\int_{\Gamma_9}g_{n1}(|w_{2n}|-|w_{1n}|)\,ds\\&\hspace{2.5cm}+\int_{\Gamma_{10}}g_{+n1}(w_{2n}-w_{1n})\,ds-\int_{\Gamma_{11}}g_{-n1}(w_{2n}-w_{1n})\,ds,\\
&J_2(w_2)-J_2(w_1)=\int_{\Gamma_8}g_{\tau
2}(|w_{2\tau}|-|w_{1\tau}|)\,ds+\int_{\Gamma_9}g_{n2}(|w_{2n}|-|w_{1n}|)\,ds\\&\hspace{2.5cm}+\int_{\Gamma_{10}}g_{+n2}(w_{2n}-w_{1n})\,ds-\int_{\Gamma_{11}}g_{-n2}(w_{2n}-w_{1n})\,ds.
\end{aligned}
\end{equation}
Subtracting two formulae in \eqref{5.48}, we have
\begin{equation}\label{5.49}
\begin{aligned}
|J_1(w_2)&-J_1(w_1)+J_2(w_1)-J_2(w_2)|\\&\leq\|g_{\tau 1}-g_{\tau
2}\|_{\textbf{L}^2_\tau(\Gamma_8)}\|w_{2\tau}-w_{1\tau}\|_{\textbf{L}^2_\tau(\Gamma_8)}+\|g_{n1}-g_{n2}\|_{L^2(\Gamma_9)}\|w_{2n}-w_{1n}\|_{L^2(\Gamma_9)}\\&
+\|g_{+n1}-g_{+n2}\|_{L^2(\Gamma_{10})}\|w_{n2}-w_{n1}\|_{L^2(\Gamma_{10})}+\|g_{-n1}-g_{-n2}\|_{L^2(\Gamma_{11})}\|w_{2n}-w_{1n}\|_{L^2(\Gamma_{11})}\\&\leq
M\big(\|g_{\tau 1}-g_{\tau
2}\|_{\textbf{L}^2_\tau(\Gamma_8)}+\|g_{n1}-g_{n2}\|_{L^2(\Gamma_9)}+\|g_{+n1}-g_{+n2}\|_{L^2(\Gamma_{10})}\\&\hspace{1cm}+\|g_{-n1}-g_{-n2}\|_{L^2(\Gamma_{11})}\big)\|w_2-w_1\|_{\mathbf{V}(\Omega)}.
\end{aligned}
\end{equation}
By \eqref{5.47}, \eqref{5.49} we have
\[
\begin{aligned}
\|w_1-w_2\|_{\mathbf{V}}\leq\f{
M}{\delta}\Big(\|F_3^1-F_3^2\|_{\mathbf{V}(\Omega)^*}&+ \|g_{\tau
1}-g_{\tau
2}\|_{\textbf{L}^2_\tau(\Gamma_8)}+\|g_{n1}-g_{n2}\|_{L^2(\Gamma_9)}\\&+\|g_{+n1}-g_{+n2}\|_{L^2(\Gamma_{10})}+\|g_{-n1}-g_{-n2}\|_{L^2(\Gamma_{11})}\Big),
\end{aligned}
\]
from which we get \eqref{5.42}. $\square$

\begin{remark}\label{r5.3}
The estimates of solutions \eqref{5.10}, \eqref{5.23}, \eqref{5.41}
are independent from thresholds $g_\tau, g_n, g_{+n},$ $ g_{-n}.$
(cf. (8) in \cite{al1}, (25) in \cite{ll3}.)
\end{remark}

{\bf Acknowledgment:}\,\,The authors are grateful to the anonymous
referee for his or her valuable comments.

\end{document}